\documentclass[12pt]{article}
\usepackage{pgf,tikz}
\usetikzlibrary{arrows}
\usepackage{amsmath}%
\usepackage{amsfonts}%
\usepackage{amssymb}%
\usepackage{graphicx}%
\newtheorem{lemma}{Lemma}

\begin{document}

\title{On Star-critical $(K_{1,n},K_{1,m}+e)$ Ramsey numbers}
\author{C. J. Jayawardene \\
Department of Mathematics\\
University of Colombo \\
Sri Lanka\\
email: c\_jayawardene@maths.cmb.ac.lk\\
\\
J. N. Senadheera,  K. A. S. N. Fernando and W. C. W. Navaratna \\
Department of Mathematics \\
The Open University of Sri Lanka \\
Sri Lanka\\
email: jnsen@ou.ac.lk, kafer@ou.ac.lk and wcper@ou.ac.lk \\
}

\maketitle
\begin{abstract}   We say that $K_n\rightarrow (G,H)$, if for every  red/blue coloring of edges of the complete graph $K_n$,  there exists a red copy of $G$, or a blue copy of $H$ in the coloring of $K_{n}$.  The Ramsey number $r(G, H)$ is the smallest positive integer $n$ such that $K_{n} \rightarrow (G, H)$.  Let $r(n,m)=r(K_n, K_m)$.  A closely related concept of Ramsey numbers is the Star-critical Ramsey number $r_*(G, H)$ defined as the largest value of $k$  such that $K_{r(G,H)-1} \sqcup K_{1,k} \rightarrow (G, H)$. Literature on survey papers in this area reveals many unsolved problems related to these numbers. One of these problems is the calculation of Ramsey numbers for certain classes of graphs. The primary objective of this paper is to calculate the Star critical Ramsey numbers for the case of Stars versus  $K_{1,m}+e.$  The methodology that we follow in solving this problem is to first find a closed form for the Ramsey number  $r_*(K_{1,n}, K_{1,m}+e)$ for all $n,m \geq 3$.  Based on the values of $r_*(K_{1,n}, K_{1,m}+e)$ for different $n, m$, we arrive at a general formula for  $r_*(K_{1,n}, K_{1,m}+e)$.  Henceforth, we show that   $r_*(K_{1,n}, K_{1,m}+e)=n+m-1$ is defined by a piecewise function related to the  three disjoint cases of $n,m$ both even and $n \leq m-2$, $n$ or $m$ is odd and $n \leq m-2$ and $n>m-2$.
\end{abstract}

\noindent Keywords: Ramsey theory, Star-critical Ramsey numbers\\
\noindent Mathematics Subject Classification: 05C55, 05C38, 05D10  \\

\section{Introduction}
\noindent Given two graphs $G$ and $H$, we say that $K_{n} \rightarrow (G, H)$, if any red and blue two colouring of $K_n$ contains a copy of $G$ (in the first color red) or a copy of $H$ (in the second color blue).   Studies on Star-critical Ramsey numbers related to different classes of graphs are trees vs complete graphs \cite{HoIs}, paths vs. paths \cite{Ho}, cycles vs. cycles \cite{ZhBrCh} and complete graphs vs stripes \cite{LiLi} are some such examples. In this paper, we extend this list by calculating Star-critical Ramsey numbers related to stars versus $K_{1,m}+e$.

\vspace{14pt}

\section{Notation}
\noindent Consider a simple graph $G$ and let $v \in V(G)$. We denote the  neighborhood of $v$ by $\Gamma(v)$ which represents the set of vertices adjacent to $v$. The degree of $v$ which is equal to $|\Gamma(v)|$ is denoted by $d(v)$.  Consider a  red/blue colouring of the complete graph $K_{n}$ given by $K_{n} =H_R \oplus H_B$ where $H_R$ and $H_B$ denote the red and blue graphs with vertex set $V(G)$. Likewise, the degree of vertex $v$ in $H_R$ and $H_B$ are denoted by $d_R(v)$ and $d_B(v)$ respectively. Then clearly, we get $n-1=d_R(v)+d_B(v)$. 

\vspace{8pt}

\section{The exact values of $r_*(K_{1,n}, K_{1,m}+e)$ for $n,m \geq 3$ }

\vspace{8pt}

\noindent  In order to find lower bounds for Star critical Ramsey numbers, we deal with constructions of graphs generated by regular $K_n$ convex $n$-gons drawn in an Euclidean plane. Label the vertices of $K_n$ by  $v_0,v_1,v_2,...,v_{n-1}$ in the anti-clockwise order. Given any $0 \leq i,k \leq n-1$, $v_{i+k\pmod{n}}$ and $v_{i-k\pmod{n}}$ are represented by the two vertices separated from $v_i$ by a path of length $k$ along the outer cycle of the $n$-gon, in the anti-clockwise direction and the clockwise direction respectively. The red/blue colorings of $K_n$ in such a scenario are called \textit{standard regular colorings of $K_n$}. The following lemma plays an crucial role in finding $r_*(K_{1,n}, K_{1,m}+e)$ for $n,m \geq 3$.

\vspace{8pt}

\begin{lemma}
\label{l1}
\noindent Given $n,m \geq 3$,

\[ r(K_{1,n},K_{1,m}+e) =
\begin{cases} 
  \hspace{4pt} n+m-1 & \text{ if } n \text{ and  } m  \text{ are both even and } n \leq m-2 \\
\hspace{20pt} & \\
  \hspace{4pt} n+m & \text{ if } n \text{ or  } m  \text{ is odd } \text{ and  } n \leq m-2 \\
\hspace{20pt} & \\
  \hspace{4pt} 2n+1 & \text{ if } n > m-2 \\
\end{cases}
\]
\end{lemma}

\vspace{10pt}
\noindent {\bf Proof.} We break up the proof in to 4 parts correspondingly.

\vspace{10pt}
\noindent  \textbf{Case 1.} \textit{If $n$ and  $m$  are both even and $n \leq m-2$ }

\vspace{8pt}
\noindent Consider a standard coloring on $K_{n+m-2}$ such that each $v_i \in V(H_R)$  ($0 \leq i \leq n+m-3$) is adjacent in red to all vertices of $\{v_{(i \pm k) \bmod (n+m-2)} $ $ \vert $ $ 0<k \leq \frac{n-2}{2}\}$  and adjacent in blue to all the other vertices  of $V(K_{n+m-2})\setminus \{v_i\}$ except for the $\frac{m-2}{2}$ diagonal red edges joining $v_i$ to the diametrically opposite vertex $v_{(i+\frac{n+m-2}{2}) \bmod (n+m-2)}$ when $i=0,1,...,\frac{m-2}{2}-1$ (see Figure 1). We note that there are many alternative colorings with different number of red diagonals. However, this particular coloring was selected as the same coloring can be used to find Star-critical Ramsey numbers. Such a  coloring is well defined, since by definition, $(v_i,v_j)$ is a red edge iff $(v_j,v_i)$ is a red edge. In such a construction, any vertex of  $K_{n+m-2}$ will  be adjacent in red  to $\frac{n-2}{2}$ vertices immediately left of it,  $\frac{n-2}{2}$ vertices immediately right of it and at most one vertex opposite it. Therefore, the red degree of any vertex adjacent in red to its opposite vertex is equal to  $2\times \frac{n-2}{2}+1=n-1$. Similarly, the red degree of any vertex not adjacent in red to its opposite vertex is equal   $2\times \frac{n-2}{2}=n-2$. Accordingly, the blue degree  will be $(n+m-3)-(n-1)=m-2$ or else $(n+m-3)-(n-2)=m-1$, respectively. In this coloring, $H_R$ has no red $K_{1,n}$. Also $H_B$ has no blue $K_{1,m}+e$. That is, $K_{n+m-2}  \not \rightarrow (K_{1,n}, K_{1,m}+e)$. Hence, $r(K_{1,n},K_{1,m}+e) \geq n+m-1$.

\vspace{8pt}
\noindent Next we need to show that, $r(K_{1,n},K_{1,m}+e) \leq n+m-1$. Suppose there exists a red/blue coloring of $K_{n+m-1}$ such that $H_R$ contains no $K_{1,n}$ and $H_B$ contains no $K_{1,m}+e$. In order to avoid a red $K_{1,n}$,  every vertex $v \in V(K_{n+m-1})$ must satisfy $d_R(v) \leq n-1$. However, by Handshaking lemma, all vertices of $V(K_{n+m-1})$ cannot have $d_R(v) = n-1$ since otherwise it will force $H_R$ to have an odd number of odd degree vertices. Therefore, there exists a vertex $v_0 \in V(K_{n+m-1})$ such that $d_R(v_0) \leq n-2$. Hence $d_B(v_0) \geq m$. In order to avoid a blue $K_{1,m}+e$, all vertices of $\Gamma_B(v_0)$ must be adjacent to each other in red. That is, the vertices of $\Gamma_B(v_0)$ induce a red complete graph of order at least $m$.

\begin{center}

\definecolor{uuuuuu}{rgb}{0.26666666666666666,0.26666666666666666,0.26666666666666666}
\begin{tikzpicture}[line cap=round,line join=round,>=triangle 45,x=0.9048223350253808cm,y=0.8910081743869209cm]
\clip(-4.699999999999998,2.7000000000000015) rectangle (11.060000000000011,10.04);
\draw (-4.819999999999998,1.8200000000000016) node[anchor=north west] {$v_{1,3}$};
\draw (-0.6399999999999958,1.6400000000000017) node[anchor=north west] {$v_{ 2,3}$};
\draw (4.547473508864644E-15,0.24000000000000218) node[anchor=north west] {$v_{2,1}$};
\draw (4.547473508864644E-15,0.24000000000000218) node[anchor=north west] {$u_{4}$};
\draw (4.547473508864644E-15,0.24000000000000218) node[anchor=north west] {$y_{2}$};
\draw (1.7764036952755078,6.296403695275507)-- (-4.256403695275508,7.496403695275508);
\draw (-1.8399999999999992,9.912807390551016)-- (0.46865543901354423,9.453587271712909);
\draw (0.46865543901354423,9.453587271712909)-- (1.7764036952755078,7.496403695275507);
\draw (1.7764036952755078,7.496403695275507)-- (1.3171835764374002,5.187748256261964);
\draw (1.3171835764374002,5.187748256261964)-- (-0.64,3.88);
\draw (-2.9486554390135455,4.339220118838108)-- (-4.256403695275508,6.296403695275507);
\draw (-4.256403695275508,6.296403695275507)-- (-3.7971835764373987,8.605059134289053);
\draw (-3.7971835764373987,8.605059134289053)-- (-1.8399999999999992,9.912807390551016);
\draw (-2.9486554390135433,9.453587271712909)-- (-0.6399999999999996,9.912807390551016);
\draw (-0.6399999999999996,9.912807390551016)-- (1.3171835764374007,8.60505913428905);
\draw (1.3171835764374007,8.60505913428905)-- (1.7764036952755078,6.296403695275507);
\draw (1.7764036952755078,6.296403695275507)-- (0.4686554390135431,4.3392201188381065);
\draw (-1.84,3.88)-- (-3.7971835764374005,5.1877482562619655);
\draw (-3.7971835764374005,5.1877482562619655)-- (-4.256403695275508,7.496403695275508);
\draw (-2.9486554390135433,9.453587271712909)-- (-4.256403695275508,7.496403695275508);
\draw [dash pattern=on 4pt off 4pt] (7.52,9.792807390551012)-- (3.903596304724494,7.376403695275506);
\draw [dash pattern=on 4pt off 4pt] (7.52,9.792807390551012)-- (3.9035963047244944,6.176403695275506);
\draw [dash pattern=on 4pt off 4pt] (7.52,9.792807390551012)-- (4.362816423562601,5.0677482562619645);
\draw [dash pattern=on 4pt off 4pt] (7.52,9.792807390551012)-- (5.211344560986456,4.219220118838108);
\draw [dash pattern=on 4pt off 4pt] (7.52,9.792807390551012)-- (7.52,3.76);
\draw [dash pattern=on 4pt off 4pt] (7.52,9.792807390551012)-- (8.628655439013542,4.219220118838107);
\draw [dash pattern=on 4pt off 4pt] (7.52,9.792807390551012)-- (9.477183576437398,5.067748256261963);
\draw [dash pattern=on 4pt off 4pt] (7.52,9.792807390551012)-- (9.936403695275505,6.176403695275506);
\draw (-1.8399999999999992,9.912807390551016)-- (-4.256403695275508,7.496403695275508);
\draw (-2.9486554390135433,9.453587271712909)-- (-4.256403695275508,6.296403695275507);
\draw (-3.7971835764373987,8.605059134289053)-- (-3.7971835764374005,5.1877482562619655);
\draw (-4.256403695275508,7.496403695275508)-- (-2.9486554390135455,4.339220118838108);
\draw (-4.256403695275508,6.296403695275507)-- (-1.84,3.88);
\draw (-2.9486554390135455,4.339220118838108)-- (0.4686554390135431,4.3392201188381065);
\draw (-1.84,3.88)-- (1.3171835764374002,5.187748256261964);
\draw (-0.64,3.88)-- (1.7764036952755078,6.296403695275507);
\draw (0.4686554390135431,4.3392201188381065)-- (1.7764036952755078,7.496403695275507);
\draw (1.3171835764374002,5.187748256261964)-- (1.3171835764374007,8.60505913428905);
\draw (1.7764036952755078,6.296403695275507)-- (0.46865543901354423,9.453587271712909);
\draw (1.7764036952755078,7.496403695275507)-- (-0.6399999999999996,9.912807390551016);
\draw (1.3171835764374007,8.60505913428905)-- (-1.8399999999999992,9.912807390551016);
\draw [dash pattern=on 4pt off 4pt] (6.320000000000001,9.792807390551012)-- (5.211344560986456,4.219220118838108);
\draw [dash pattern=on 4pt off 4pt] (6.320000000000001,9.792807390551012)-- (4.362816423562601,5.0677482562619645);
\draw [dash pattern=on 4pt off 4pt] (6.320000000000001,9.792807390551012)-- (3.9035963047244944,6.176403695275506);
\draw [dash pattern=on 4pt off 4pt] (6.320000000000001,9.792807390551012)-- (6.32,3.76);
\draw [dash pattern=on 4pt off 4pt] (6.320000000000001,9.792807390551012)-- (8.628655439013542,4.219220118838107);
\draw [dash pattern=on 4pt off 4pt] (6.320000000000001,9.792807390551012)-- (9.477183576437398,5.067748256261963);
\draw [dash pattern=on 4pt off 4pt] (6.320000000000001,9.792807390551012)-- (9.936403695275505,6.176403695275506);
\draw [dash pattern=on 4pt off 4pt] (6.320000000000001,9.792807390551012)-- (9.936403695275505,7.3764036952755045);
\draw [dash pattern=on 4pt off 4pt] (4.362816423562601,5.0677482562619645)-- (5.211344560986458,9.333587271712904);
\draw [dash pattern=on 4pt off 4pt] (5.211344560986458,9.333587271712904)-- (5.211344560986456,4.219220118838108);
\draw [dash pattern=on 4pt off 4pt] (5.211344560986458,9.333587271712904)-- (6.32,3.76);
\draw [dash pattern=on 4pt off 4pt] (5.211344560986458,9.333587271712904)-- (7.52,3.76);
\draw [dash pattern=on 4pt off 4pt] (5.211344560986458,9.333587271712904)-- (9.4771835764374,8.485059134289047);
\draw [dash pattern=on 4pt off 4pt] (5.211344560986458,9.333587271712904)-- (9.936403695275505,7.3764036952755045);
\draw [dash pattern=on 4pt off 4pt] (5.211344560986458,9.333587271712904)-- (9.936403695275505,6.176403695275506);
\draw [dash pattern=on 4pt off 4pt] (5.211344560986458,9.333587271712904)-- (9.477183576437398,5.067748256261963);
\draw [dash pattern=on 4pt off 4pt] (4.362816423562602,8.48505913428905)-- (8.628655439013542,9.333587271712904);
\draw [dash pattern=on 4pt off 4pt] (4.362816423562602,8.48505913428905)-- (9.4771835764374,8.485059134289047);
\draw [dash pattern=on 4pt off 4pt] (4.362816423562602,8.48505913428905)-- (9.936403695275505,7.3764036952755045);
\draw [dash pattern=on 4pt off 4pt] (4.362816423562602,8.48505913428905)-- (9.936403695275505,6.176403695275506);
\draw [dash pattern=on 4pt off 4pt] (4.362816423562602,8.48505913428905)-- (8.628655439013542,4.219220118838107);
\draw [dash pattern=on 4pt off 4pt] (4.362816423562602,8.48505913428905)-- (7.52,3.76);
\draw [dash pattern=on 4pt off 4pt] (4.362816423562602,8.48505913428905)-- (6.32,3.76);
\draw [dash pattern=on 4pt off 4pt] (4.362816423562602,8.48505913428905)-- (5.211344560986456,4.219220118838108);
\draw [dash pattern=on 4pt off 4pt] (3.903596304724494,7.376403695275506)-- (8.628655439013542,9.333587271712904);
\draw [dash pattern=on 4pt off 4pt] (3.903596304724494,7.376403695275506)-- (9.4771835764374,8.485059134289047);
\draw [dash pattern=on 4pt off 4pt] (3.903596304724494,7.376403695275506)-- (9.936403695275505,7.3764036952755045);
\draw [dash pattern=on 4pt off 4pt] (3.903596304724494,7.376403695275506)-- (9.477183576437398,5.067748256261963);
\draw [dash pattern=on 4pt off 4pt] (3.903596304724494,7.376403695275506)-- (8.628655439013542,4.219220118838107);
\draw [dash pattern=on 4pt off 4pt] (3.903596304724494,7.376403695275506)-- (7.52,3.76);
\draw [dash pattern=on 4pt off 4pt] (3.903596304724494,7.376403695275506)-- (6.32,3.76);
\draw [dash pattern=on 4pt off 4pt] (3.9035963047244944,6.176403695275506)-- (8.628655439013542,9.333587271712904);
\draw [dash pattern=on 4pt off 4pt] (3.9035963047244944,6.176403695275506)-- (9.4771835764374,8.485059134289047);
\draw [dash pattern=on 4pt off 4pt] (3.9035963047244944,6.176403695275506)-- (9.936403695275505,6.176403695275506);
\draw [dash pattern=on 4pt off 4pt] (3.9035963047244944,6.176403695275506)-- (9.477183576437398,5.067748256261963);
\draw [dash pattern=on 4pt off 4pt] (3.9035963047244944,6.176403695275506)-- (8.628655439013542,4.219220118838107);
\draw [dash pattern=on 4pt off 4pt] (3.9035963047244944,6.176403695275506)-- (7.52,3.76);
\draw [dash pattern=on 4pt off 4pt] (4.362816423562601,5.0677482562619645)-- (8.628655439013542,9.333587271712904);
\draw [dash pattern=on 4pt off 4pt] (4.362816423562601,5.0677482562619645)-- (9.936403695275505,7.3764036952755045);
\draw [dash pattern=on 4pt off 4pt] (4.362816423562601,5.0677482562619645)-- (9.936403695275505,6.176403695275506);
\draw [dash pattern=on 4pt off 4pt] (4.362816423562601,5.0677482562619645)-- (9.477183576437398,5.067748256261963);
\draw [dash pattern=on 4pt off 4pt] (4.362816423562601,5.0677482562619645)-- (8.628655439013542,4.219220118838107);
\draw [dash pattern=on 4pt off 4pt] (5.211344560986456,4.219220118838108)-- (9.4771835764374,8.485059134289047);
\draw [dash pattern=on 4pt off 4pt] (5.211344560986456,4.219220118838108)-- (9.936403695275505,7.3764036952755045);
\draw [dash pattern=on 4pt off 4pt] (5.211344560986456,4.219220118838108)-- (9.936403695275505,6.176403695275506);
\draw [dash pattern=on 4pt off 4pt] (5.211344560986456,4.219220118838108)-- (9.477183576437398,5.067748256261963);
\draw [dash pattern=on 4pt off 4pt] (6.32,3.76)-- (8.628655439013542,9.333587271712904);
\draw [dash pattern=on 4pt off 4pt] (6.32,3.76)-- (9.4771835764374,8.485059134289047);
\draw [dash pattern=on 4pt off 4pt] (6.32,3.76)-- (9.936403695275505,7.3764036952755045);
\draw [dash pattern=on 4pt off 4pt] (6.32,3.76)-- (9.936403695275505,6.176403695275506);
\draw [dash pattern=on 4pt off 4pt] (7.52,3.76)-- (8.628655439013542,9.333587271712904);
\draw [dash pattern=on 4pt off 4pt] (7.52,3.76)-- (9.4771835764374,8.485059134289047);
\draw [dash pattern=on 4pt off 4pt] (7.52,3.76)-- (9.936403695275505,7.3764036952755045);
\draw [dash pattern=on 4pt off 4pt] (8.628655439013542,4.219220118838107)-- (8.628655439013542,9.333587271712904);
\draw [dash pattern=on 4pt off 4pt] (8.628655439013542,4.219220118838107)-- (9.4771835764374,8.485059134289047);
\draw [dash pattern=on 4pt off 4pt] (9.477183576437398,5.067748256261963)-- (8.628655439013542,9.333587271712904);
\draw (-1.5599999999999963,3.780000000000001) node[anchor=north west] {$H_R$};
\draw (6.620000000000008,3.760000000000001) node[anchor=north west] {$H_B$};
\draw (-0.64,3.88)-- (-3.7971835764374005,5.1877482562619655);
\draw (-0.64,3.88)-- (-2.9486554390135455,4.339220118838108);
\draw (-2.9486554390135433,9.453587271712909)-- (-1.8399999999999992,9.912807390551016);
\draw (-1.8399999999999992,9.912807390551016)-- (-0.6399999999999996,9.912807390551016);
\draw (-0.6399999999999996,9.912807390551016)-- (0.46865543901354423,9.453587271712909);
\draw (0.46865543901354423,9.453587271712909)-- (1.3171835764374007,8.60505913428905);
\draw (1.3171835764374007,8.60505913428905)-- (1.7764036952755078,7.496403695275507);
\draw (1.7764036952755078,7.496403695275507)-- (1.7764036952755078,6.296403695275507);
\draw (1.7764036952755078,6.296403695275507)-- (1.3171835764374002,5.187748256261964);
\draw (1.3171835764374002,5.187748256261964)-- (0.4686554390135431,4.3392201188381065);
\draw (0.4686554390135431,4.3392201188381065)-- (-0.64,3.88);
\draw (-0.64,3.88)-- (-1.84,3.88);
\draw (-1.84,3.88)-- (-2.9486554390135455,4.339220118838108);
\draw (-2.9486554390135455,4.339220118838108)-- (-3.7971835764374005,5.1877482562619655);
\draw (-3.7971835764374005,5.1877482562619655)-- (-4.256403695275508,6.296403695275507);
\draw (-4.256403695275508,6.296403695275507)-- (-4.256403695275508,7.496403695275508);
\draw (-4.256403695275508,7.496403695275508)-- (-3.7971835764373987,8.605059134289053);
\draw (-3.7971835764373987,8.605059134289053)-- (-2.9486554390135433,9.453587271712909);
\draw [dash pattern=on 4pt off 4pt] (6.320000000000001,9.792807390551012)-- (7.52,3.76);
\draw [dash pattern=on 4pt off 4pt] (9.477183576437398,5.067748256261963)-- (4.362816423562602,8.48505913428905);
\draw (2.040000000000006,6.72) node[anchor=north west] {$v_0$};
\draw (10.20000000000001,6.68) node[anchor=north west] {$v_0$};
\draw (1.4400000000000053,9.319999999999999) node[anchor=north west] {$v_2$};
\draw (9.68000000000001,9.12) node[anchor=north west] {$v_2$};
\draw (1.6200000000000054,5.460000000000001) node[anchor=north west] {$v_{15}$};
\draw (2.040000000000006,8.08) node[anchor=north west] {$v_1$};
\draw (10.18000000000001,7.9399999999999995) node[anchor=north west] {$v_1$};
\draw (10.00000000000001,5.380000000000001) node[anchor=north west] {$v_{15}$};
\draw (1.7764036952755078,7.496403695275507)-- (-4.256403695275508,6.296403695275507);
\draw (1.3171835764374007,8.60505913428905)-- (-3.7971835764374005,5.1877482562619655);
\draw [dash pattern=on 4pt off 4pt] (7.52,9.792807390551012)-- (6.32,3.76);
\draw [dash pattern=on 4pt off 4pt] (5.211344560986458,9.333587271712904)-- (8.628655439013542,4.219220118838107);
\draw (0.46865543901354423,9.453587271712909)-- (-2.9486554390135455,4.339220118838108);
\draw (-3.7971835764373987,8.605059134289053)-- (-0.6399999999999996,9.912807390551016);
\draw (-2.9486554390135433,9.453587271712909)-- (0.46865543901354423,9.453587271712909);
\begin{scriptsize}
\draw [fill=black] (-1.84,3.88) circle (1.5pt);
\draw [fill=black] (-0.64,3.88) circle (1.5pt);
\draw [fill=uuuuuu] (0.4686554390135431,4.3392201188381065) circle (1.5pt);
\draw [fill=uuuuuu] (1.3171835764374002,5.187748256261964) circle (1.5pt);
\draw [fill=uuuuuu] (1.7764036952755078,6.296403695275507) circle (1.5pt);
\draw [fill=uuuuuu] (1.7764036952755078,7.496403695275507) circle (1.5pt);
\draw [fill=uuuuuu] (1.3171835764374007,8.60505913428905) circle (1.5pt);
\draw [fill=uuuuuu] (0.46865543901354423,9.453587271712909) circle (1.5pt);
\draw [fill=uuuuuu] (-0.6399999999999996,9.912807390551016) circle (1.5pt);
\draw [fill=uuuuuu] (-1.8399999999999992,9.912807390551016) circle (1.5pt);
\draw [fill=uuuuuu] (-2.9486554390135433,9.453587271712909) circle (1.5pt);
\draw [fill=uuuuuu] (-3.7971835764373987,8.605059134289053) circle (1.5pt);
\draw [fill=uuuuuu] (-4.256403695275508,7.496403695275508) circle (1.5pt);
\draw [fill=uuuuuu] (-4.256403695275508,6.296403695275507) circle (1.5pt);
\draw [fill=uuuuuu] (-3.7971835764374005,5.1877482562619655) circle (1.5pt);
\draw [fill=uuuuuu] (-2.9486554390135455,4.339220118838108) circle (1.5pt);
\draw [fill=black] (6.32,3.76) circle (1.5pt);
\draw [fill=black] (7.52,3.76) circle (1.5pt);
\draw [fill=uuuuuu] (8.628655439013542,4.219220118838107) circle (1.5pt);
\draw [fill=uuuuuu] (9.477183576437398,5.067748256261963) circle (1.5pt);
\draw [fill=uuuuuu] (9.936403695275505,6.176403695275506) circle (1.5pt);
\draw [fill=uuuuuu] (9.936403695275505,7.3764036952755045) circle (1.5pt);
\draw [fill=uuuuuu] (9.4771835764374,8.485059134289047) circle (1.5pt);
\draw [fill=uuuuuu] (8.628655439013542,9.333587271712904) circle (1.5pt);
\draw [fill=uuuuuu] (7.52,9.792807390551012) circle (1.5pt);
\draw [fill=uuuuuu] (6.320000000000001,9.792807390551012) circle (1.5pt);
\draw [fill=uuuuuu] (5.211344560986458,9.333587271712904) circle (1.5pt);
\draw [fill=uuuuuu] (4.362816423562602,8.48505913428905) circle (1.5pt);
\draw [fill=uuuuuu] (3.903596304724494,7.376403695275506) circle (1.5pt);
\draw [fill=uuuuuu] (3.9035963047244944,6.176403695275506) circle (1.5pt);
\draw [fill=uuuuuu] (4.362816423562601,5.0677482562619645) circle (1.5pt);
\draw [fill=uuuuuu] (5.211344560986456,4.219220118838108) circle (1.5pt);
\end{scriptsize}
\end{tikzpicture}
\end{center}
\begin{center}
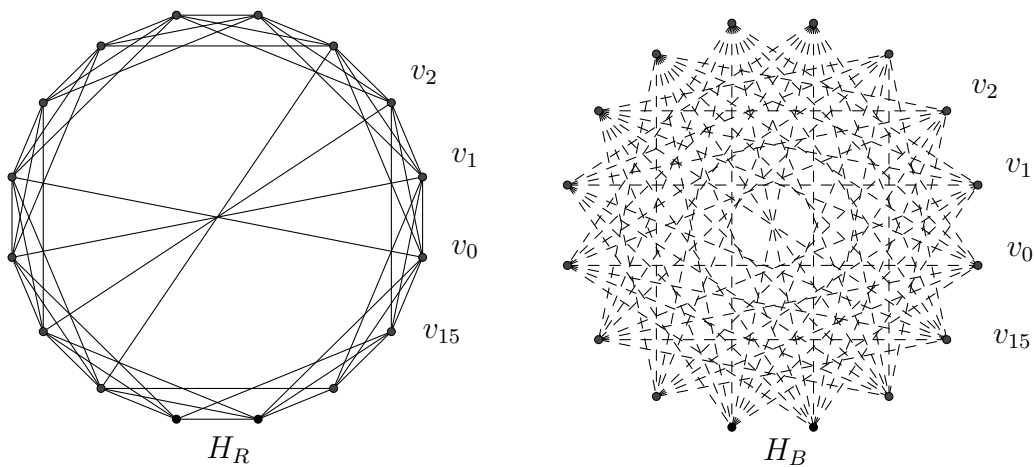
Figure 1. A Ramsey critical $(K_{1,8}, K_{1,10}+e)$ coloring of $K_{16} =H_R \oplus H_B$
\end{center}
\vspace{8pt}

\noindent  Let $w \in \Gamma_B(v_0)$. Then, $d_R(w) \geq m-1 \geq n$. That is, $H_R$ contains a red $K_{1,n}$, a contradiction. Therefore, $K_{n+m-1}  \rightarrow (K_{1,n}, K_{1,m}+e)$. Hence, $r(K_{1,n},K_{1,m}+e) \leq n+m-1$.   Combining with the earlier result, we find $r(K_{1,n},K_{1,m}+e) = n+m-1$, as required.

\vspace{26pt}
\noindent  \textbf{Case 2.} \textit{If $n$ is odd and $n \leq m-2$ }

\vspace{8pt}

\noindent As before, consider a standard coloring on $K_{n+m-1}$ such that each $v_i \in V(H_R)$  ($0 \leq i \leq n+m-2)$ is adjacent to $\{v_{(i \pm k) \bmod (n+m-1)} $ $ \vert $ $ 0<k \leq \frac{n-1}{2}\}$ in red and adjacent to all the other vertices  of $V(K_{n+m-1})\setminus \{v_i\}$ in blue. This coloring is also  well defined. In such a construction, any vertex of  $K_{n+m-1}$ will  be adjacent in red to $\frac{n-1}{2}$ vertices immediately left of it, $\frac{n-1}{2}$ vertices immediately right of it.  The red degree of any vertex is equal to  $2\times \frac{n-1}{2}=n-1$ and the blue degree of any vertex is $(n+m-2)-n-1=m-1$. Therefore, $H_R$ has no red $K_{1,n}$. Also $H_B$ has no blue $K_{1,m}+e$. That is, $K_{n+m-1}  \not \rightarrow (K_{1,n}, K_{1,m}+e)$. Hence, $r(K_{1,n},K_{1,m}+e) \geq n+m$.

\vspace{8pt}
\noindent Next we need to show that, $r(K_{1,n},K_{1,m}+e) \leq n+m$. Suppose there exists a red/blue coloring of $K_{n+m}$ such that $H_R$ contains no $K_{1,n}$ and $H_B$ contains no $K_{1,m}+e$. In order to avoid a red $K_{1,n}$,  every vertex $v \in V(K_{n+m})$ must satisfy $d_R(v) \leq n-1$. That is, for any vertex  $v\in V(K_{n+m})$,  $d_B(v) \geq m$. Let $v_0 \in V(K_{n+m})$. In order to avoid a blue $K_{1,m}+e$, all vertices of $\Gamma_B(v_0)$ must be adjacent to each other in red. However, as $n+1\leq m$, we argue that $\Gamma_B(v_0)$ contains a red $K_{1,n}$, a contradiction. Hence, $r(K_{1,n},K_{1,m}+e) \leq n+m$. Combining with the earlier result, $r(K_{1,n},K_{1,m}+e) = n+m$, as required.

\vspace{16pt}
\noindent  \textbf{Case 3.} \textit{If $n$ is even, $m$ is odd and $n \leq m-2$ }

\vspace{8pt}
\noindent Now consider a standard coloring on $K_{n+m-1}$ such that each $v_i \in V(H_B)$  ($0 \leq i \leq n+m-2)$ is adjacent to $\{v_{(i \pm k) \bmod (n+m-1)} $ $ \vert $ $ 0<k \leq \frac{m-1}{2}\}$ in blue and adjacent to all the other vertices of $V(K_{n+m-1})\setminus \{v_i\}$ in red. This coloring is also  well defined. In such a construction, any vertex of  $K_{n+m-1}$ will  be adjacent in blue to $\frac{m-1}{2}$ vertices immediately left of it, $\frac{m-1}{2}$ vertices immediately right of it. Therefore, the blue degree of any vertex is equal to  $2\times \frac{m-1}{2}=m-1$ and the red degree of any vertex is $(n+m-2)-(m-1)=n-1$. Therefore, $H_R$ has no red $K_{1,n}$. Also $H_B$ has no $K_{1,m}+e$ since it has no blue $K_{1,m}$. That is, $K_{n+m-1}  \not \rightarrow (K_{1,n}, K_{1,m}+e)$. Hence, $r(K_{1,n},K_{1,m}+e) \geq n+m$.

\vspace{8pt}
\noindent Next we need to show that, $r(K_{1,n},K_{1,m}+e) \leq n+m$. Suppose there exists a red/blue coloring of $K_{n+m}$ such that $H_R$ contains no $K_{1,n}$ and $H_B$ contains no $K_{1,m}+e$. In order to avoid a red $K_{1,n}$,  every vertex $v \in V(K_{n+m})$ must satisfy $d_R(v)\leq n-1$. Hence, for any vertex  $v\in V(K_{n+m})$,  $d_B(v) \geq m$. Let $v_0 \in V(K_{n+m})$. In order to avoid a blue $K_{1,m}+e$, all vertices of $\Gamma_B(v_0)$ must be adjacent to each other in red. As $n+1\leq m$, $\Gamma_B(v_0)$ contains a red $K_{1,n}$, a contradiction. Hence, $r(K_{1,n},K_{1,m}+e) = n+m$.

\vspace{16pt}
\noindent  \textbf{Case 4.} \textit{$n > m-2$ }

\vspace{8pt}
\noindent Consider a standard regular coloring of $K_{2n}=H_R \oplus H_B$ such that each $v_i \in V(K_{2n})$  ($0 \leq i \leq n)$ forms a red clique of size $n$ and each $v_i \in V(K_{2n})$  ($n+1 \leq i \leq 2n)$ also forms an independent red clique of size $n$. That is, $H_R\cong 2K_n$ and $H_B \cong K_{n,n}$ (see Figure 2).  

\begin{center}

\begin{tikzpicture}[line cap=round,line join=round,>=triangle 45,x=1.0cm,y=1.0cm]
\clip(-5.576363636363642,0.6218181818181893) rectangle (-0.7945454545454534,5.785454545454557);
\draw [dash pattern=on 5pt off 5pt] (-5.19454545454546,5.003636363636375)-- (-5.230909090909094,1.5672727272727371);
\draw [dash pattern=on 5pt off 5pt] (-5.19454545454546,5.003636363636375)-- (-4.394545454545457,1.5309090909091005);
\draw [dash pattern=on 5pt off 5pt] (-5.19454545454546,5.003636363636375)-- (-3.5763636363636375,1.5127272727272825);
\draw [dash pattern=on 5pt off 5pt] (-4.430909090909093,4.985454545454559)-- (-5.230909090909094,1.5672727272727371);
\draw [dash pattern=on 5pt off 5pt] (-4.430909090909093,4.985454545454559)-- (-4.394545454545457,1.5309090909091005);
\draw [dash pattern=on 5pt off 5pt] (-4.430909090909093,4.985454545454559)-- (-3.5763636363636375,1.5127272727272825);
\draw [dash pattern=on 5pt off 5pt] (-3.667272727272729,4.985454545454559)-- (-5.230909090909094,1.5672727272727371);
\draw [dash pattern=on 5pt off 5pt] (-3.667272727272729,4.985454545454559)-- (-4.394545454545457,1.5309090909091005);
\draw [dash pattern=on 5pt off 5pt] (-3.667272727272729,4.985454545454559)-- (-3.5763636363636375,1.5127272727272825);
\draw [dash pattern=on 5pt off 5pt] (-5.19454545454546,5.003636363636375)-- (-1.2672727272727256,1.476363636363646);
\draw [dash pattern=on 5pt off 5pt] (-4.430909090909093,4.985454545454559)-- (-1.2672727272727256,1.476363636363646);
\draw [dash pattern=on 5pt off 5pt] (-3.667272727272729,4.985454545454559)-- (-1.2672727272727256,1.476363636363646);
\draw [dash pattern=on 5pt off 5pt] (-1.2672727272727256,1.476363636363646)-- (-1.2854545454545439,5.003636363636377);
\draw [dash pattern=on 5pt off 5pt] (-1.2672727272727256,1.476363636363646)-- (-2.9218181818181823,5.003636363636377);
\draw [dash pattern=on 5pt off 5pt] (-1.2672727272727256,1.476363636363646)-- (-2.085454545454545,5.003636363636377);
\draw [dash pattern=on 5pt off 5pt] (-2.085454545454545,5.003636363636377)-- (-2.758181818181822,1.4945454545454546);
\draw [dash pattern=on 5pt off 5pt] (-2.9218181818181823,5.003636363636377)-- (-2.758181818181822,1.4945454545454546);
\draw [dash pattern=on 5pt off 5pt] (-2.9218181818181823,5.003636363636377)-- (-1.94,1.4945454545454622);
\draw [dash pattern=on 5pt off 5pt] (-2.085454545454545,5.003636363636377)-- (-1.94,1.4945454545454622);
\draw [dash pattern=on 5pt off 5pt] (-2.9218181818181823,5.003636363636377)-- (-5.230909090909094,1.5672727272727371);
\draw [dash pattern=on 5pt off 5pt] (-2.085454545454545,5.003636363636377)-- (-5.230909090909094,1.5672727272727371);
\draw [dash pattern=on 5pt off 5pt] (-1.94,1.4945454545454622)-- (-1.2854545454545439,5.003636363636377);
\draw [dash pattern=on 5pt off 5pt] (-1.94,1.4945454545454622)-- (-3.667272727272729,4.985454545454559);
\draw [dash pattern=on 5pt off 5pt] (-1.94,1.4945454545454622)-- (-4.430909090909093,4.985454545454559);
\draw [dash pattern=on 5pt off 5pt] (-1.94,1.4945454545454622)-- (-5.19454545454546,5.003636363636375);
\draw [dash pattern=on 5pt off 5pt] (-2.758181818181822,1.4945454545454546)-- (-1.2854545454545439,5.003636363636377);
\draw [dash pattern=on 5pt off 5pt] (-2.758181818181822,1.4945454545454546)-- (-3.667272727272729,4.985454545454559);
\draw [dash pattern=on 5pt off 5pt] (-2.758181818181822,1.4945454545454546)-- (-4.430909090909093,4.985454545454559);
\draw [dash pattern=on 5pt off 5pt] (-2.758181818181822,1.4945454545454546)-- (-5.19454545454546,5.003636363636375);
\draw [dash pattern=on 5pt off 5pt] (-3.5763636363636375,1.5127272727272825)-- (-2.9218181818181823,5.003636363636377);
\draw [dash pattern=on 5pt off 5pt] (-3.5763636363636375,1.5127272727272825)-- (-2.085454545454545,5.003636363636377);
\draw [dash pattern=on 5pt off 5pt] (-3.5763636363636375,1.5127272727272825)-- (-1.2854545454545439,5.003636363636377);
\draw [dash pattern=on 5pt off 5pt] (-4.394545454545457,1.5309090909091005)-- (-2.9218181818181823,5.003636363636377);
\draw [dash pattern=on 5pt off 5pt] (-4.394545454545457,1.5309090909091005)-- (-2.085454545454545,5.003636363636377);
\draw [dash pattern=on 5pt off 5pt] (-4.394545454545457,1.5309090909091005)-- (-1.2854545454545439,5.003636363636377);
\draw [dash pattern=on 5pt off 5pt] (-5.230909090909094,1.5672727272727371)-- (-1.2854545454545439,5.003636363636377);
\draw (-5.19454545454546,5.003636363636375)-- (-4.430909090909093,4.985454545454559);
\draw (-4.430909090909093,4.985454545454559)-- (-1.2854545454545439,5.003636363636377);
\draw [shift={(-4.440683131897008,4.173526011560707)}] plot[domain=0.809689847466171:2.3080944070070037,variable=\t]({1.0*1.1213347395546562*cos(\t r)+-0.0*1.1213347395546562*sin(\t r)},{0.0*1.1213347395546562*cos(\t r)+1.0*1.1213347395546562*sin(\t r)});
\draw [shift={(-2.8664153322405754,4.129042995839124)}] plot[domain=0.8418946175809634:2.3226855294388913,variable=\t]({1.0*1.1725243321536643*cos(\t r)+-0.0*1.1725243321536643*sin(\t r)},{0.0*1.1725243321536643*cos(\t r)+1.0*1.1725243321536643*sin(\t r)});
\draw [shift={(-4.0581818181818194,3.6712252964427012)}] plot[domain=0.8646426914377155:2.27694996215208,variable=\t]({1.0*1.7511829047902965*cos(\t r)+-0.0*1.7511829047902965*sin(\t r)},{0.0*1.7511829047902965*cos(\t r)+1.0*1.7511829047902965*sin(\t r)});
\draw [shift={(-3.6400000000000015,3.289862258953178)}] plot[domain=0.834078486342882:2.3075141672469126,variable=\t]({1.0*2.3137919638832356*cos(\t r)+-0.0*2.3137919638832356*sin(\t r)},{0.0*2.3137919638832356*cos(\t r)+1.0*2.3137919638832356*sin(\t r)});
\draw [shift={(-3.2399999999999998,2.685454545454556)}] plot[domain=0.8702999568471705:2.2712926967426244,variable=\t]({1.0*3.0321963782105343*cos(\t r)+-0.0*3.0321963782105343*sin(\t r)},{0.0*3.0321963782105343*cos(\t r)+1.0*3.0321963782105343*sin(\t r)});
\draw [shift={(-2.103636363636363,3.9133333333333478)}] plot[domain=0.9270284564927336:2.2145641970970598,variable=\t]({1.0*1.363151563653608*cos(\t r)+-0.0*1.363151563653608*sin(\t r)},{0.0*1.363151563653608*cos(\t r)+1.0*1.363151563653608*sin(\t r)});
\draw [shift={(-2.4665229615745083,3.705417057169648)}] plot[domain=0.8326148768520374:2.3242446557727416,variable=\t]({1.0*1.7550771980853368*cos(\t r)+-0.0*1.7550771980853368*sin(\t r)},{0.0*1.7550771980853368*cos(\t r)+1.0*1.7550771980853368*sin(\t r)});
\draw [shift={(-2.843124090423984,2.389558552440098)}] plot[domain=1.0334137613897836:2.1197394570873334,variable=\t]({1.0*3.04298163226695*cos(\t r)+-0.0*3.04298163226695*sin(\t r)},{0.0*3.04298163226695*cos(\t r)+1.0*3.04298163226695*sin(\t r)});
\draw [shift={(-3.655254237288137,3.242465331278902)}] plot[domain=1.1761915888027186:1.989496284495447,variable=\t]({1.0*1.907787161262603*cos(\t r)+-0.0*1.907787161262603*sin(\t r)},{0.0*1.907787161262603*cos(\t r)+1.0*1.907787161262603*sin(\t r)});
\draw [shift={(-3.248947745168219,3.803350035790994)}] plot[domain=0.8009622871697555:2.356133931844451,variable=\t]({1.0*1.6716470005932818*cos(\t r)+-0.0*1.6716470005932818*sin(\t r)},{0.0*1.6716470005932818*cos(\t r)+1.0*1.6716470005932818*sin(\t r)});
\draw (-5.230909090909094,1.5672727272727371)-- (-1.2672727272727256,1.476363636363646);
\draw [shift={(-3.201505308525176,3.596550366484154)}] plot[domain=3.9269597371464364:5.4519547050299915,variable=\t]({1.0*2.8699211572720276*cos(\t r)+-0.0*2.8699211572720276*sin(\t r)},{0.0*2.8699211572720276*cos(\t r)+1.0*2.8699211572720276*sin(\t r)});
\draw [shift={(-3.9415989684074746,3.33108961960029)}] plot[domain=4.081173965348956:5.284797419012405,variable=\t]({1.0*2.1848044813312626*cos(\t r)+-0.0*2.1848044813312626*sin(\t r)},{0.0*2.1848044813312626*cos(\t r)+1.0*2.1848044813312626*sin(\t r)});
\draw [shift={(-2.0015006615006676,2.4060366660367167)}] plot[domain=4.019527966854847:5.380860959018042,variable=\t]({1.0*1.18464450427212*cos(\t r)+-0.0*1.18464450427212*sin(\t r)},{0.0*1.18464450427212*cos(\t r)+1.0*1.18464450427212*sin(\t r)});
\draw [shift={(-2.783243054097558,4.2364891408309555)}] plot[domain=4.175249545635541:5.214648230982434,variable=\t]({1.0*3.149040938480208*cos(\t r)+-0.0*3.149040938480208*sin(\t r)},{0.0*3.149040938480208*cos(\t r)+1.0*3.149040938480208*sin(\t r)});
\draw [shift={(-3.5325783893525733,3.4830633882247453)}] plot[domain=4.296587112646589:5.083753717469341,variable=\t]({1.0*2.133985384674687*cos(\t r)+-0.0*2.133985384674687*sin(\t r)},{0.0*2.133985384674687*cos(\t r)+1.0*2.133985384674687*sin(\t r)});
\draw [shift={(-3.1333591231463576,3.8018955512571866)}] plot[domain=4.205449584213515:5.18970091433206,variable=\t]({1.0*2.5976855979106186*cos(\t r)+-0.0*2.5976855979106186*sin(\t r)},{0.0*2.5976855979106186*cos(\t r)+1.0*2.5976855979106186*sin(\t r)});
\draw [shift={(-4.357758799584965,2.931619442892501)}] plot[domain=4.143092224018955:5.215775541521828,variable=\t]({1.0*1.6198251114431605*cos(\t r)+-0.0*1.6198251114431605*sin(\t r)},{0.0*1.6198251114431605*cos(\t r)+1.0*1.6198251114431605*sin(\t r)});
\draw [shift={(-3.529885549455663,4.045406159858526)}] plot[domain=4.110839983430501:5.269746275571187,variable=\t]({1.0*3.0057655261643745*cos(\t r)+-0.0*3.0057655261643745*sin(\t r)},{0.0*3.0057655261643745*cos(\t r)+1.0*3.0057655261643745*sin(\t r)});
\draw [shift={(-2.740095693779903,3.131387559808631)}] plot[domain=4.235516555731295:5.167040097242866,variable=\t]({1.0*1.8219234881760282*cos(\t r)+-0.0*1.8219234881760282*sin(\t r)},{0.0*1.8219234881760282*cos(\t r)+1.0*1.8219234881760282*sin(\t r)});
\draw [shift={(-2.382426302453494,3.9959297942030925)}] plot[domain=4.2642145772385955:5.129069923831046,variable=\t]({1.0*2.755318696624454*cos(\t r)+-0.0*2.755318696624454*sin(\t r)},{0.0*2.755318696624454*cos(\t r)+1.0*2.755318696624454*sin(\t r)});
\begin{scriptsize}
\draw [fill=black] (-5.19454545454546,5.003636363636375) circle (1.5pt);
\draw [fill=black] (-4.430909090909093,4.985454545454559) circle (1.5pt);
\draw [fill=black] (-3.667272727272729,4.985454545454559) circle (1.5pt);
\draw [fill=black] (-2.085454545454545,5.003636363636377) circle (1.5pt);
\draw [fill=black] (-1.2854545454545439,5.003636363636377) circle (1.5pt);
\draw [fill=black] (-5.230909090909094,1.5672727272727371) circle (1.5pt);
\draw [fill=black] (-4.394545454545457,1.5309090909091005) circle (1.5pt);
\draw [fill=black] (-3.5763636363636375,1.5127272727272825) circle (1.5pt);
\draw [fill=black] (-1.94,1.4945454545454622) circle (1.5pt);
\draw [fill=black] (-1.2672727272727256,1.476363636363646) circle (1.5pt);
\draw [fill=black] (-2.9218181818181823,5.003636363636377) circle (1.5pt);
\draw [fill=black] (-2.758181818181822,1.4945454545454546) circle (1.5pt);
\end{scriptsize}
\end{tikzpicture}

\end{center}
\begin{center}
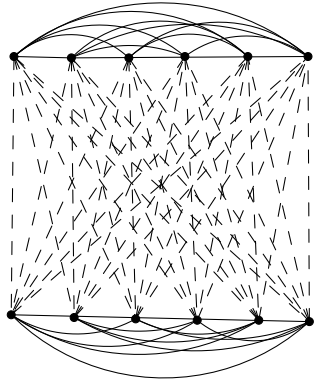
Figure 2. A Red/blue graph corresponding to a coloring of $K_{12}$ with no red $K_{1,6}$ and no blue $K_{1,7}+e$ 
\end{center}

\noindent Clearly,  $H_R$ has no $K_{1,n}$. Furthermore, $H_B$ has no $K_{1,m}+e$, since it has no blue $C_3$. That is, $K_{2n}  \not \rightarrow (K_{1,n}, K_{1,m}+e)$. Hence, $r(K_{1,n},K_{1,m}+e) \geq 2n+1$.

\vspace{8pt}
\noindent Next we need to show that, $r(K_{1,n},K_{1,m}+e) \leq 2n+1$. Suppose there exists a red/blue coloring of $K_{2n+1}$ such that $H_R$ contains no $K_{1,n}$ and $H_B$ contains no $K_{1,m}+e$. 

\begin{center}

\begin{tikzpicture}[line cap=round,line join=round,>=triangle 45,x=1.0cm,y=1.0cm]
\clip(-5.067272727272728,0.385454545454556) rectangle (8.496363636363656,7.530909090909109);
\draw (-4.321818181818185,0.9309090909091)-- (-4.176363636363639,4.149090909090922);
\draw (-4.176363636363639,4.149090909090922)-- (-3.812727272727275,0.9854545454545546);
\draw (-4.176363636363639,4.149090909090922)-- (-1.9763636363636377,1.7127272727272826);
\draw (-4.176363636363639,4.149090909090922)-- (-1.54,2.04);
\draw (-4.176363636363639,4.149090909090922)-- (-4.7945454545454576,1.0036363636363728);
\draw [dash pattern=on 4pt off 4pt] (-4.176363636363639,4.149090909090922)-- (-0.9581818181818187,3.4218181818181934);
\draw [dash pattern=on 4pt off 4pt] (-4.176363636363639,4.149090909090922)-- (-0.54,4.04);
\draw [dash pattern=on 4pt off 4pt] (-4.176363636363639,4.149090909090922)-- (-0.612727272727273,4.6581818181818315);
\draw [dash pattern=on 4pt off 4pt] (-4.176363636363639,4.149090909090922)-- (-0.8854545454545458,5.312727272727287);
\draw [dash pattern=on 4pt off 4pt] (-4.176363636363639,4.149090909090922)-- (-2.6490909090909107,6.8581818181818335);
\draw [dash pattern=on 4pt off 4pt] (-4.176363636363639,4.149090909090922)-- (-2.04909090909091,6.603636363636379);
\draw (-0.9581818181818187,3.4218181818181934)-- (-0.54,4.04);
\draw (-0.9581818181818187,3.4218181818181934)-- (-0.612727272727273,4.6581818181818315);
\draw (-0.9581818181818187,3.4218181818181934)-- (-0.8854545454545458,5.312727272727287);
\draw (-0.9581818181818187,3.4218181818181934)-- (-2.04909090909091,6.603636363636379);
\draw (-0.9581818181818187,3.4218181818181934)-- (-2.6490909090909107,6.8581818181818335);
\draw (-5.194545454545456,4.821818181818197) node[anchor=north west] {$v_0$};
\draw [rotate around={-59.713635558172236:(-1.603636363636365,4.9490909090909225)}] (-1.603636363636365,4.9490909090909225) ellipse (2.7168078923503547cm and 1.0502123358556898cm);
\draw (-0.6127272727272672,7.167272727272746) node[anchor=north west] {Blue neigbourhood will be forced to induce a red };
\draw (-0.6490909090909036,6.749090909090927) node[anchor=north west] {$K_{n+1}$};
\begin{scriptsize}
\draw [fill=black] (-4.176363636363639,4.149090909090922) circle (1.5pt);
\draw [fill=black] (-4.321818181818185,0.9309090909091) circle (1.5pt);
\draw [fill=black] (-3.812727272727275,0.9854545454545546) circle (1.5pt);
\draw [fill=black] (-1.9763636363636377,1.7127272727272826) circle (1.5pt);
\draw [fill=black] (-1.54,2.04) circle (1.5pt);
\draw [fill=black] (-4.7945454545454576,1.0036363636363728) circle (1.5pt);
\draw [fill=black] (-0.9581818181818187,3.4218181818181934) circle (1.5pt);
\draw [fill=black] (-0.54,4.04) circle (1.5pt);
\draw [fill=black] (-0.612727272727273,4.6581818181818315) circle (1.5pt);
\draw [fill=black] (-0.8854545454545458,5.312727272727287) circle (1.5pt);
\draw [fill=black] (-2.6490909090909107,6.8581818181818335) circle (1.5pt);
\draw [fill=black] (-2.04909090909091,6.603636363636379) circle (1.5pt);
\draw [fill=black] (-3.230909090909093,1.1127272727272821) circle (1.5pt);
\draw [fill=black] (-2.794545454545456,1.276363636363646) circle (1.5pt);
\draw [fill=black] (-2.3581818181818197,1.4945454545454642) circle (1.5pt);
\draw [fill=black] (-1.3218181818181824,6.003636363636378) circle (1.5pt);
\draw [fill=black] (-1.0854545454545461,5.64) circle (1.5pt);
\draw [fill=black] (-1.6127272727272737,6.294545454545469) circle (1.5pt);
\end{scriptsize}
\end{tikzpicture}

\end{center}
\begin{center}
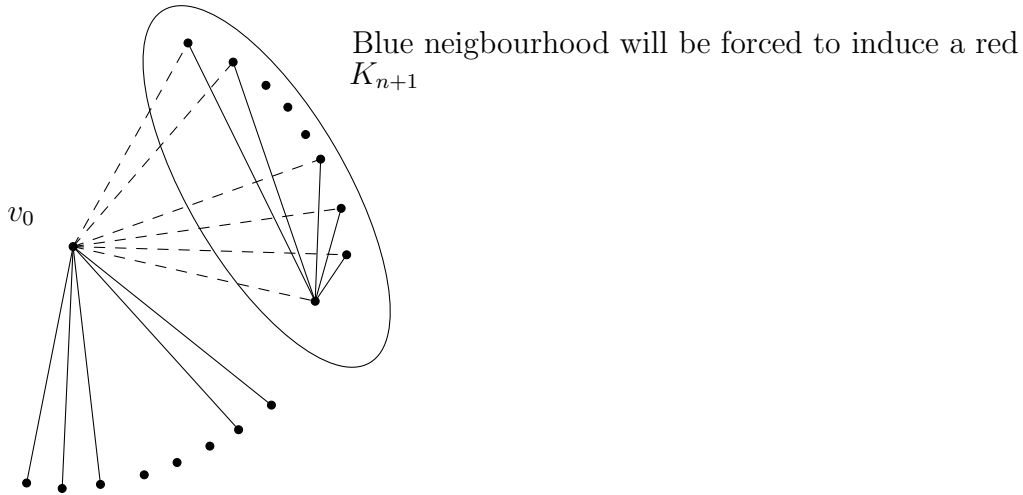
Figure 3.  Neighborhood of a vertex of $K_{2n+1}$ used in the argument containing no red $K_{1,n}$
\end{center} 

\vspace{8pt}
\noindent Let $v_0 \in V(K_{2n+1})$. In order to avoid a red $K_{1,n}$,  $v_0$ must satisfy $d_B(v_0) \geq 2n-(n-1) = n+1 \geq m$. To avoid a blue $K_{1,m}+e$, all vertices of $\Gamma_B(v_0)$ must be adjacent to each other in red. That is, the vertices of $\Gamma_B(v_0)$ induces a red complete graph of order at least $n+1$ (see Figure 3). Hence, $V(K_{2n+1})$ will contain a vertex of red degree $n$, a contradiction.

\begin{lemma}
\label{l2}
\noindent Given $n,m \geq 3$,

\[ r_*(K_{1,n},K_{1,m}+e) =
\begin{cases} 
  \hspace{4pt} n+m-2 & \text{ if } n \text{ and  } m  \text{ are both even } n \leq m-2 \\
\hspace{20pt} & \\
  \hspace{4pt} 1 & \text{ if } n \text{ or  } m  \text{ is odd } \text{ and  } n \leq m-2 \\
\hspace{20pt} & \\
  \hspace{4pt} n+1 & \text{ if } n > m-2 \\
\end{cases}
\]
\end{lemma}

\vspace{10pt}
\noindent {\bf Proof.} We break up the proof in to 3 cases.

\vspace{10pt}
\noindent  \textbf{Case 1.} \textit{$n$ and $m$ are even and $n \leq m-2$ }

\vspace{8pt}
\noindent To show that, $r_*(K_{1,n},K_{1,m}+e)\geq n+m-2$, consider the coloring of $K_{n+m-2} \sqcup K_{1,n+m-3}$ introduced in  Case 1 of Lemma 1. Add a vertex (say $x$) and connect it in blue to the vertices $v_i$ and the diametrically opposite vertices $v_{j \bmod (n+m-2)}$ for $i=0,1,...,\frac{m-2}{2}-1$ where $j=i+\frac{n+m-2}{2}$. Connect all the other vertices excluding $v_{\frac{n+m-4}{2}}$ to $x$ in red (see Figure 4). 

\begin{center}

\begin{tikzpicture}[line cap=round,line join=round,>=triangle 45,x=0.7593120881534575cm,y=0.7359169781469995cm]
\clip(-2.264557864248511,-3.2957727653620235) rectangle (10.199715796128483,7.306953530153469);
\draw (9.066600008821483,2.612616697024472) node[anchor=north west] {$x$};
\draw (0.84,-0.76)-- (0.02,-0.16);
\draw (-0.5209368419504998,4.077189474148501)-- (0.14,5.126978947316833);
\draw (0.14,5.126978947316833)-- (1.1,5.76);
\draw (1.1,5.76)-- (2.18,5.98);
\draw (2.18,5.98)-- (3.2,5.94);
\draw (3.2,5.94)-- (4.28,5.573021052683167);
\draw (4.28,5.573021052683167)-- (5.14,4.9);
\draw (0.84,-0.76)-- (2.1,-1.1069789473168332);
\draw (2.1,-1.1069789473168332)-- (3.16,-1.14);
\draw (3.16,-1.14)-- (4.2,-0.78);
\draw (5.18,-0.07302105268316642)-- (4.2,-0.78);
\draw (5.18,-0.07302105268316642)-- (5.86,0.86);
\draw [line width=1.2000000000000002pt] (3.16,-1.14)-- (2.18,5.98);
\draw [line width=1.2000000000000002pt] (0.84,-0.76)-- (4.28,5.573021052683167);
\draw [dash pattern=on 5pt off 5pt] (5.18,-0.07302105268316642)-- (0.14,5.126978947316833);
\draw [line width=1.2000000000000002pt] (3.2,5.94)-- (2.1,-1.1069789473168332);
\draw [line width=1.2000000000000002pt] (5.14,4.9)-- (0.02,-0.16);
\draw [dash pattern=on 5pt off 5pt] (1.1,5.76)-- (4.2,-0.78);
\draw [dash pattern=on 5pt off 5pt] (5.86,0.86)-- (-0.5209368419504998,4.077189474148501);
\draw (2.18,5.98)-- (0.14,5.126978947316833);
\draw (1.1,5.76)-- (-0.5209368419504998,4.077189474148501);
\draw (1.1,5.76)-- (3.2,5.94);
\draw (1.1,5.76)-- (4.28,5.573021052683167);
\draw (3.2,5.94)-- (5.14,4.9);
\draw (-0.5209368419504998,4.077189474148501)-- (-0.88,3.26);
\draw (-0.88,3.26)-- (-1.0609368419504999,2.0860421053663334);
\draw (-1.0609368419504999,2.0860421053663334)-- (-0.78,0.88);
\draw (-0.78,0.88)-- (0.02,-0.16);
\draw (5.14,4.9)-- (5.94,4.02);
\draw (5.94,4.02)-- (6.326978947316833,2.98);
\draw (6.326978947316833,2.98)-- (6.273021052683167,1.96);
\draw (6.273021052683167,1.96)-- (5.86,0.86);
\draw [dash pattern=on 5pt off 5pt] (-0.88,3.26)-- (6.273021052683167,1.96);
\draw [line width=1.2000000000000002pt] (-0.78,0.88)-- (5.94,4.02);
\draw [line width=1.2000000000000002pt] (-1.0609368419504999,2.0860421053663334)-- (6.326978947316833,2.98);
\draw (0.14,5.126978947316833)-- (-0.88,3.26);
\draw (-0.5209368419504998,4.077189474148501)-- (-1.0609368419504999,2.0860421053663334);
\draw (-0.88,3.26)-- (-0.78,0.88);
\draw (-1.0609368419504999,2.0860421053663334)-- (0.02,-0.16);
\draw (-0.78,0.88)-- (0.84,-0.76);
\draw (0.02,-0.16)-- (2.1,-1.1069789473168332);
\draw (0.84,-0.76)-- (3.16,-1.14);
\draw (2.1,-1.1069789473168332)-- (4.2,-0.78);
\draw (3.16,-1.14)-- (5.18,-0.07302105268316642);
\draw (4.2,-0.78)-- (5.86,0.86);
\draw (6.273021052683167,1.96)-- (5.18,-0.07302105268316642);
\draw (5.86,0.86)-- (6.326978947316833,2.98);
\draw (5.94,4.02)-- (6.273021052683167,1.96);
\draw (5.14,4.9)-- (6.326978947316833,2.98);
\draw (2.18,5.98)-- (5.14,4.9);
\draw (3.2,5.94)-- (5.94,4.02);
\draw (4.28,5.573021052683167)-- (6.326978947316833,2.98);
\draw (5.14,4.9)-- (6.273021052683167,1.96);
\draw (5.94,4.02)-- (5.86,0.86);
\draw (6.326978947316833,2.98)-- (5.18,-0.07302105268316642);
\draw (6.273021052683167,1.96)-- (4.2,-0.78);
\draw (5.86,0.86)-- (3.16,-1.14);
\draw (5.18,-0.07302105268316642)-- (2.1,-1.1069789473168332);
\draw (4.2,-0.78)-- (0.84,-0.76);
\draw (3.16,-1.14)-- (0.02,-0.16);
\draw (2.1,-1.1069789473168332)-- (-0.78,0.88);
\draw (0.84,-0.76)-- (-1.0609368419504999,2.0860421053663334);
\draw (0.02,-0.16)-- (-0.88,3.26);
\draw (-0.78,0.88)-- (-0.5209368419504998,4.077189474148501);
\draw (-1.0609368419504999,2.0860421053663334)-- (0.14,5.126978947316833);
\draw (-0.88,3.26)-- (1.1,5.76);
\draw (-0.5209368419504998,4.077189474148501)-- (2.18,5.98);
\draw [dash pattern=on 5pt off 5pt] (0.14,5.126978947316833)-- (3.2,5.94);
\draw [dash pattern=on 5pt off 5pt] (8.634936851752151,2.153974592638306)-- (6.326978947316833,2.98);
\draw [dash pattern=on 5pt off 5pt] (8.634936851752151,2.153974592638306)-- (5.14,4.9);
\draw [shift={(4.029495128860938,1.3367530683291158)},dash pattern=on 5pt off 5pt]  plot[domain=0.1756189016381606:1.7490808842597274,variable=\t]({1.0*4.677386501321163*cos(\t r)+-0.0*4.677386501321163*sin(\t r)},{0.0*4.677386501321163*cos(\t r)+1.0*4.677386501321163*sin(\t r)});
\draw [shift={(3.584246193587072,1.2756422155765617)}] plot[domain=0.17218147321631802:2.0766903483180634,variable=\t]({1.0*5.126494307913638*cos(\t r)+-0.0*5.126494307913638*sin(\t r)},{0.0*5.126494307913638*cos(\t r)+1.0*5.126494307913638*sin(\t r)});
\draw [shift={(3.8930998662634724,4.502771925307103)}] plot[domain=4.770418283096239:5.823276740754068,variable=\t]({1.0*5.2916790253094925*cos(\t r)+-0.0*5.2916790253094925*sin(\t r)},{0.0*5.2916790253094925*cos(\t r)+1.0*5.2916790253094925*sin(\t r)});
\draw [shift={(3.8301118649267365,3.60435375063848)},dash pattern=on 5pt off 5pt]  plot[domain=4.360453473220939:5.9900239787833955,variable=\t]({1.0*5.018958343718658*cos(\t r)+-0.0*5.018958343718658*sin(\t r)},{0.0*5.018958343718658*cos(\t r)+1.0*5.018958343718658*sin(\t r)});
\draw [shift={(3.868012160058487,3.022799163220892)},dash pattern=on 5pt off 5pt]  plot[domain=4.037362337959161:6.1029031706161065,variable=\t]({1.0*4.845454277023647*cos(\t r)+-0.0*4.845454277023647*sin(\t r)},{0.0*4.845454277023647*cos(\t r)+1.0*4.845454277023647*sin(\t r)});
\draw [shift={(3.8731766423326874,1.9182148534334833)},dash pattern=on 5pt off 5pt]  plot[domain=-2.9220687339832265:0.04947065138885134,variable=\t]({1.0*4.767592992972494*cos(\t r)+-0.0*4.767592992972494*sin(\t r)},{0.0*4.767592992972494*cos(\t r)+1.0*4.767592992972494*sin(\t r)});
\draw [shift={(3.3802722282565285,-0.056014613829853795)},dash pattern=on 5pt off 5pt]  plot[domain=0.3981180663925302:1.4123000696570709,variable=\t]({1.0*5.700487049202214*cos(\t r)+-0.0*5.700487049202214*sin(\t r)},{0.0*5.700487049202214*cos(\t r)+1.0*5.700487049202214*sin(\t r)});
\draw [shift={(3.7058845964134974,1.1962232532968955)},dash pattern=on 5pt off 5pt]  plot[domain=0.1919159642818591:1.8795652984178184,variable=\t]({1.0*5.02123926574599*cos(\t r)+-0.0*5.02123926574599*sin(\t r)},{0.0*5.02123926574599*cos(\t r)+1.0*5.02123926574599*sin(\t r)});
\draw [shift={(3.8980653598488364,2.242077156329902)}] plot[domain=-0.018597170455498357:2.4868870718970157,variable=\t]({1.0*4.7376907447125935*cos(\t r)+-0.0*4.7376907447125935*sin(\t r)},{0.0*4.7376907447125935*cos(\t r)+1.0*4.7376907447125935*sin(\t r)});
\draw [shift={(3.9429155127782716,2.4286793611476063)},dash pattern=on 5pt off 5pt]  plot[domain=3.7248866617880068:6.224704853411036,variable=\t]({1.0*4.70005605873249*cos(\t r)+-0.0*4.70005605873249*sin(\t r)},{0.0*4.70005605873249*cos(\t r)+1.0*4.70005605873249*sin(\t r)});
\draw [dash pattern=on 5pt off 5pt] (8.634936851752151,2.153974592638306)-- (5.94,4.02);
\draw (8.634936851752151,2.153974592638306)-- (5.18,-0.07302105268316642);
\draw [dash pattern=on 5pt off 5pt] (1.1,5.76)-- (-1.0609368419504999,2.0860421053663334);
\draw [dash pattern=on 5pt off 5pt] (1.1,5.76)-- (-0.78,0.88);
\draw [dash pattern=on 5pt off 5pt] (1.1,5.76)-- (0.02,-0.16);
\draw [dash pattern=on 5pt off 5pt] (1.1,5.76)-- (0.84,-0.76);
\draw [dash pattern=on 5pt off 5pt] (1.1,5.76)-- (2.1,-1.1069789473168332);
\draw [dash pattern=on 5pt off 5pt] (1.1,5.76)-- (3.16,-1.14);
\draw [dash pattern=on 5pt off 5pt] (1.1,5.76)-- (5.18,-0.07302105268316642);
\draw [dash pattern=on 5pt off 5pt] (1.1,5.76)-- (5.86,0.86);
\draw [dash pattern=on 5pt off 5pt] (1.1,5.76)-- (6.273021052683167,1.96);
\draw [dash pattern=on 5pt off 5pt] (1.1,5.76)-- (6.326978947316833,2.98);
\draw [dash pattern=on 5pt off 5pt] (1.1,5.76)-- (5.94,4.02);
\draw [dash pattern=on 5pt off 5pt] (1.1,5.76)-- (5.14,4.9);
\draw [dash pattern=on 5pt off 5pt] (2.18,5.98)-- (-0.88,3.26);
\draw [dash pattern=on 5pt off 5pt] (2.18,5.98)-- (-1.0609368419504999,2.0860421053663334);
\draw [dash pattern=on 5pt off 5pt] (2.18,5.98)-- (-0.78,0.88);
\draw [dash pattern=on 5pt off 5pt] (2.18,5.98)-- (0.02,-0.16);
\draw [dash pattern=on 5pt off 5pt] (2.18,5.98)-- (0.84,-0.76);
\draw [dash pattern=on 5pt off 5pt] (2.18,5.98)-- (2.1,-1.1069789473168332);
\draw [dash pattern=on 5pt off 5pt] (2.18,5.98)-- (4.2,-0.78);
\draw [dash pattern=on 5pt off 5pt] (2.18,5.98)-- (5.18,-0.07302105268316642);
\draw [dash pattern=on 5pt off 5pt] (2.18,5.98)-- (5.86,0.86);
\draw [dash pattern=on 5pt off 5pt] (2.18,5.98)-- (6.273021052683167,1.96);
\draw [dash pattern=on 5pt off 5pt] (2.18,5.98)-- (6.326978947316833,2.98);
\draw [dash pattern=on 5pt off 5pt] (2.18,5.98)-- (5.94,4.02);
\draw [dash pattern=on 5pt off 5pt] (3.2,5.94)-- (-0.88,3.26);
\draw [dash pattern=on 5pt off 5pt] (3.2,5.94)-- (-0.5209368419504998,4.077189474148501);
\draw [dash pattern=on 5pt off 5pt] (3.2,5.94)-- (-1.0609368419504999,2.0860421053663334);
\draw [dash pattern=on 5pt off 5pt] (3.2,5.94)-- (-0.78,0.88);
\draw [dash pattern=on 5pt off 5pt] (3.2,5.94)-- (0.02,-0.16);
\draw [dash pattern=on 5pt off 5pt] (3.2,5.94)-- (0.84,-0.76);
\draw [dash pattern=on 5pt off 5pt] (3.2,5.94)-- (3.16,-1.14);
\draw [dash pattern=on 5pt off 5pt] (3.2,5.94)-- (4.2,-0.78);
\draw [dash pattern=on 5pt off 5pt] (3.2,5.94)-- (5.18,-0.07302105268316642);
\draw [dash pattern=on 5pt off 5pt] (3.2,5.94)-- (5.86,0.86);
\draw [dash pattern=on 5pt off 5pt] (3.2,5.94)-- (6.273021052683167,1.96);
\draw [dash pattern=on 5pt off 5pt] (3.2,5.94)-- (6.326978947316833,2.98);
\draw [dash pattern=on 5pt off 5pt] (4.28,5.573021052683167)-- (6.273021052683167,1.96);
\draw [dash pattern=on 5pt off 5pt] (4.28,5.573021052683167)-- (5.86,0.86);
\draw [dash pattern=on 5pt off 5pt] (4.28,5.573021052683167)-- (5.18,-0.07302105268316642);
\draw [dash pattern=on 5pt off 5pt] (4.28,5.573021052683167)-- (4.2,-0.78);
\draw [dash pattern=on 5pt off 5pt] (4.28,5.573021052683167)-- (3.16,-1.14);
\draw [dash pattern=on 5pt off 5pt] (4.28,5.573021052683167)-- (2.1,-1.1069789473168332);
\draw [dash pattern=on 5pt off 5pt] (4.28,5.573021052683167)-- (0.02,-0.16);
\draw [dash pattern=on 5pt off 5pt] (4.28,5.573021052683167)-- (-0.78,0.88);
\draw [dash pattern=on 5pt off 5pt] (4.28,5.573021052683167)-- (-1.0609368419504999,2.0860421053663334);
\draw [dash pattern=on 5pt off 5pt] (4.28,5.573021052683167)-- (-0.88,3.26);
\draw [dash pattern=on 5pt off 5pt] (4.28,5.573021052683167)-- (-0.5209368419504998,4.077189474148501);
\draw [dash pattern=on 5pt off 5pt] (4.28,5.573021052683167)-- (0.14,5.126978947316833);
\draw [dash pattern=on 5pt off 5pt] (5.14,4.9)-- (0.14,5.126978947316833);
\draw [dash pattern=on 5pt off 5pt] (5.14,4.9)-- (-0.5209368419504998,4.077189474148501);
\draw [dash pattern=on 5pt off 5pt] (5.14,4.9)-- (-0.88,3.26);
\draw [dash pattern=on 5pt off 5pt] (5.14,4.9)-- (-1.0609368419504999,2.0860421053663334);
\draw [dash pattern=on 5pt off 5pt] (5.14,4.9)-- (-0.78,0.88);
\draw [dash pattern=on 5pt off 5pt] (5.14,4.9)-- (0.84,-0.76);
\draw [dash pattern=on 5pt off 5pt] (5.14,4.9)-- (2.1,-1.1069789473168332);
\draw [dash pattern=on 5pt off 5pt] (5.14,4.9)-- (3.16,-1.14);
\draw [dash pattern=on 5pt off 5pt] (5.14,4.9)-- (4.2,-0.78);
\draw [dash pattern=on 5pt off 5pt] (5.14,4.9)-- (5.18,-0.07302105268316642);
\draw [dash pattern=on 5pt off 5pt] (5.14,4.9)-- (5.86,0.86);
\draw [dash pattern=on 5pt off 5pt] (5.94,4.02)-- (0.14,5.126978947316833);
\draw [dash pattern=on 5pt off 5pt] (5.94,4.02)-- (-0.5209368419504998,4.077189474148501);
\draw [dash pattern=on 5pt off 5pt] (5.94,4.02)-- (-0.88,3.26);
\draw [dash pattern=on 5pt off 5pt] (5.94,4.02)-- (-1.0609368419504999,2.0860421053663334);
\draw [dash pattern=on 5pt off 5pt] (5.94,4.02)-- (0.02,-0.16);
\draw [dash pattern=on 5pt off 5pt] (5.94,4.02)-- (0.84,-0.76);
\draw [dash pattern=on 5pt off 5pt] (5.94,4.02)-- (2.1,-1.1069789473168332);
\draw [dash pattern=on 5pt off 5pt] (5.94,4.02)-- (3.16,-1.14);
\draw [dash pattern=on 5pt off 5pt] (5.94,4.02)-- (4.2,-0.78);
\draw [dash pattern=on 5pt off 5pt] (5.94,4.02)-- (5.18,-0.07302105268316642);
\draw [dash pattern=on 5pt off 5pt] (6.326978947316833,2.98)-- (0.14,5.126978947316833);
\draw [dash pattern=on 5pt off 5pt] (6.326978947316833,2.98)-- (-0.5209368419504998,4.077189474148501);
\draw [dash pattern=on 5pt off 5pt] (6.326978947316833,2.98)-- (-0.88,3.26);
\draw [dash pattern=on 5pt off 5pt] (6.326978947316833,2.98)-- (-0.78,0.88);
\draw [dash pattern=on 5pt off 5pt] (6.326978947316833,2.98)-- (0.02,-0.16);
\draw [dash pattern=on 5pt off 5pt] (6.326978947316833,2.98)-- (0.84,-0.76);
\draw [dash pattern=on 5pt off 5pt] (6.326978947316833,2.98)-- (2.1,-1.1069789473168332);
\draw [dash pattern=on 5pt off 5pt] (6.326978947316833,2.98)-- (3.16,-1.14);
\draw [dash pattern=on 5pt off 5pt] (6.326978947316833,2.98)-- (4.2,-0.78);
\draw [dash pattern=on 5pt off 5pt] (6.273021052683167,1.96)-- (0.14,5.126978947316833);
\draw [dash pattern=on 5pt off 5pt] (6.273021052683167,1.96)-- (-0.5209368419504998,4.077189474148501);
\draw [dash pattern=on 5pt off 5pt] (6.273021052683167,1.96)-- (-1.0609368419504999,2.0860421053663334);
\draw [dash pattern=on 5pt off 5pt] (6.273021052683167,1.96)-- (-0.78,0.88);
\draw [dash pattern=on 5pt off 5pt] (6.273021052683167,1.96)-- (0.02,-0.16);
\draw [dash pattern=on 5pt off 5pt] (6.273021052683167,1.96)-- (0.84,-0.76);
\draw [dash pattern=on 5pt off 5pt] (6.273021052683167,1.96)-- (2.1,-1.1069789473168332);
\draw [dash pattern=on 5pt off 5pt] (6.273021052683167,1.96)-- (3.16,-1.14);
\draw [dash pattern=on 5pt off 5pt] (5.86,0.86)-- (-0.88,3.26);
\draw [dash pattern=on 5pt off 5pt] (5.86,0.86)-- (-1.0609368419504999,2.0860421053663334);
\draw [dash pattern=on 5pt off 5pt] (5.86,0.86)-- (-0.78,0.88);
\draw [dash pattern=on 5pt off 5pt] (5.86,0.86)-- (0.02,-0.16);
\draw [dash pattern=on 5pt off 5pt] (5.86,0.86)-- (0.84,-0.76);
\draw [dash pattern=on 5pt off 5pt] (5.86,0.86)-- (2.1,-1.1069789473168332);
\draw [dash pattern=on 5pt off 5pt] (5.18,-0.07302105268316642)-- (-0.5209368419504998,4.077189474148501);
\draw [dash pattern=on 5pt off 5pt] (5.18,-0.07302105268316642)-- (-0.88,3.26);
\draw [dash pattern=on 5pt off 5pt] (5.18,-0.07302105268316642)-- (-1.0609368419504999,2.0860421053663334);
\draw [dash pattern=on 5pt off 5pt] (5.18,-0.07302105268316642)-- (-0.78,0.88);
\draw [dash pattern=on 5pt off 5pt] (5.18,-0.07302105268316642)-- (0.02,-0.16);
\draw [dash pattern=on 5pt off 5pt] (5.18,-0.07302105268316642)-- (0.84,-0.76);
\draw [dash pattern=on 5pt off 5pt] (4.2,-0.78)-- (0.14,5.126978947316833);
\draw [dash pattern=on 5pt off 5pt] (4.2,-0.78)-- (-0.5209368419504998,4.077189474148501);
\draw [dash pattern=on 5pt off 5pt] (4.2,-0.78)-- (-0.88,3.26);
\draw [dash pattern=on 5pt off 5pt] (4.2,-0.78)-- (-1.0609368419504999,2.0860421053663334);
\draw [dash pattern=on 5pt off 5pt] (4.2,-0.78)-- (-0.78,0.88);
\draw [dash pattern=on 5pt off 5pt] (4.2,-0.78)-- (0.02,-0.16);
\draw [dash pattern=on 5pt off 5pt] (3.16,-1.14)-- (0.14,5.126978947316833);
\draw [dash pattern=on 5pt off 5pt] (3.16,-1.14)-- (-0.5209368419504998,4.077189474148501);
\draw [dash pattern=on 5pt off 5pt] (3.16,-1.14)-- (-0.88,3.26);
\draw [dash pattern=on 5pt off 5pt] (3.16,-1.14)-- (-1.0609368419504999,2.0860421053663334);
\draw [dash pattern=on 5pt off 5pt] (3.16,-1.14)-- (-0.78,0.88);
\draw [dash pattern=on 5pt off 5pt] (2.1,-1.1069789473168332)-- (0.14,5.126978947316833);
\draw [dash pattern=on 5pt off 5pt] (2.1,-1.1069789473168332)-- (-0.5209368419504998,4.077189474148501);
\draw [dash pattern=on 5pt off 5pt] (2.1,-1.1069789473168332)-- (-0.88,3.26);
\draw [dash pattern=on 5pt off 5pt] (2.1,-1.1069789473168332)-- (-1.0609368419504999,2.0860421053663334);
\draw [dash pattern=on 5pt off 5pt] (0.84,-0.76)-- (0.14,5.126978947316833);
\draw [dash pattern=on 5pt off 5pt] (0.84,-0.76)-- (-0.5209368419504998,4.077189474148501);
\draw [dash pattern=on 5pt off 5pt] (0.84,-0.76)-- (-0.88,3.26);
\draw [dash pattern=on 5pt off 5pt] (0.02,-0.16)-- (0.14,5.126978947316833);
\draw [dash pattern=on 5pt off 5pt] (0.02,-0.16)-- (-0.5209368419504998,4.077189474148501);
\draw [dash pattern=on 5pt off 5pt] (-0.78,0.88)-- (0.14,5.126978947316833);
\draw [shift={(3.79026605881983,1.653850733593549)},dash pattern=on 5pt off 5pt]  plot[domain=-3.230447597174398:0.1028673735892064,variable=\t]({1.0*4.870416713822064*cos(\t r)+-0.0*4.870416713822064*sin(\t r)},{0.0*4.870416713822064*cos(\t r)+1.0*4.870416713822064*sin(\t r)});
\draw [shift={(3.9281360284627653,2.502097651966055)}] plot[domain=-0.07382729044099001:2.8013348404277294,variable=\t]({1.0*4.719657154344302*cos(\t r)+-0.0*4.719657154344302*sin(\t r)},{0.0*4.719657154344302*cos(\t r)+1.0*4.719657154344302*sin(\t r)});
\draw (6.206831593237151,2.0460588033709723) node[anchor=north west] {$v_{19}$};
\draw (-1.6979999705950113,3.988543010182971) node[anchor=north west] {$v_{9}$};
\draw (6.341726329821318,3.448964063846305) node[anchor=north west] {$v_{0}$};
\draw (5.86,0.86)-- (8.634936851752151,2.153974592638306);
\draw [shift={(8.299474826998173,-8.238124231072627)}] plot[domain=1.5385270438762528:1.7669497335136484,variable=\t]({1.0*10.397511852930304*cos(\t r)+-0.0*10.397511852930304*sin(\t r)},{0.0*10.397511852930304*cos(\t r)+1.0*10.397511852930304*sin(\t r)});
\draw [shift={(3.55820407816829,4.395093965957214)},dash pattern=on 5pt off 5pt]  plot[domain=4.64057100885667:5.867464980290036,variable=\t]({1.0*5.549399219721519*cos(\t r)+-0.0*5.549399219721519*sin(\t r)},{0.0*5.549399219721519*cos(\t r)+1.0*5.549399219721519*sin(\t r)});
\begin{scriptsize}
\draw [fill=black] (8.634936851752151,2.153974592638306) circle (1.5pt);
\draw [fill=black] (2.1,-1.1069789473168332) circle (1.5pt);
\draw [fill=black] (3.16,-1.14) circle (1.5pt);
\draw [fill=black] (4.2,-0.78) circle (1.5pt);
\draw [fill=black] (5.18,-0.07302105268316642) circle (1.5pt);
\draw [fill=black] (5.86,0.86) circle (1.5pt);
\draw [fill=black] (6.273021052683167,1.96) circle (1.5pt);
\draw [fill=black] (6.326978947316833,2.98) circle (1.5pt);
\draw [fill=black] (0.84,-0.76) circle (1.5pt);
\draw [fill=black] (0.02,-0.16) circle (1.5pt);
\draw [fill=black] (-0.5209368419504998,4.077189474148501) circle (1.5pt);
\draw [fill=black] (0.14,5.126978947316833) circle (1.5pt);
\draw [fill=black] (1.1,5.76) circle (1.5pt);
\draw [fill=black] (2.18,5.98) circle (1.5pt);
\draw [fill=black] (3.2,5.94) circle (1.5pt);
\draw [fill=black] (4.28,5.573021052683167) circle (1.5pt);
\draw [fill=black] (5.14,4.9) circle (1.5pt);
\draw [fill=black] (-0.78,0.88) circle (1.5pt);
\draw [fill=black] (-1.0609368419504999,2.0860421053663334) circle (1.5pt);
\draw [fill=black] (-0.88,3.26) circle (1.5pt);
\draw [fill=black] (5.94,4.02) circle (1.5pt);
\end{scriptsize}
\end{tikzpicture} 
\end{center}
\begin{center}
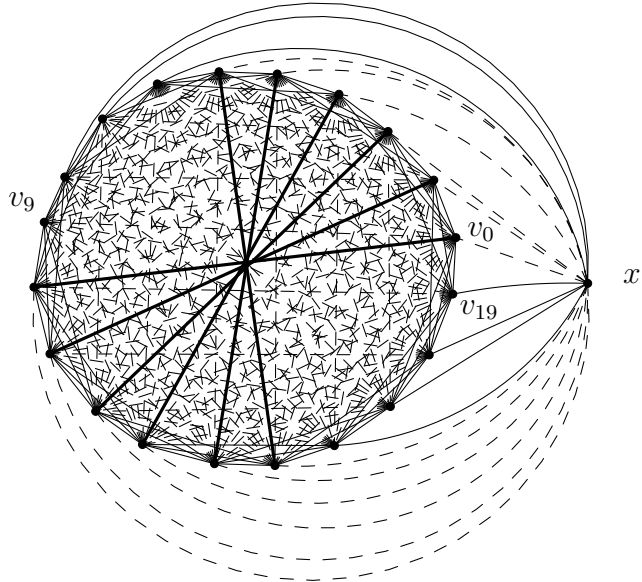
Figure 4. A red/blue  coloring of $K_{n+m-2} \sqcup K_{1,n+m-3}$ when $n=8$ and $m=14$ 
\end{center}

\noindent This coloring of $K_{n+m-2} \sqcup K_{1,n+m-3}$ contains neither red $K_{1,n}$ nor blue $K_{1,m}+e$. Thus, $K_{n+m-2} \sqcup K_{1,n+m-3}\not \rightarrow (K_{1,n}, K_{1,m}+e)$ . Therefore, $r_*(K_{1,n},K_{1,m}+e)\geq n+m-2$. Finally, using $r_*(K_{1,n},K_{1,m}+e)\leq r(K_{1,n},K_{1,m}+e)-1=n+m-2$, we  conclude that $r_*(K_{1,n},K_{1,m}+e)= n+m-2$.

\vspace{14pt}
\noindent  \textbf{Case 2.} \textit{$n$ or $m$ is odd and $n \leq m-2$ }

\vspace{8pt}
\noindent We first show that, $r_*(K_{1,n},K_{1,m}+e) \leq 1$. Suppose there exists a red/blue coloring of $K_{n+m-1} \sqcup K_{1,1}$ such that $H_R$ contains no $K_{1,n}$ and $H_B$ contains no $K_{1,m}+e$. First let us restrict our attention to the red/blue coloring of $K_{n+m-1}$. In order to avoid a red $K_{1,n}$, any vertex $v \in K_{n+m-1}$ must satisfy $d_R(v) \leq n-1$ and hence $d_B(v) \geq m-1$. Suppose that there exists a vertex $v_0 \in K_{n+m-1}$ such that $d_R(v_0) \leq n-2$. That is,  $d_B(v_0) \geq m$. In order to avoid a blue $K_{1,m}+e$, $\Gamma_B(v_0)$ must induce a red complete graph. Since $n<m-1$, the vertices of $\Gamma_B(v_0)$  will contain a red complete graph of order at least $n+1$. Hence, $H_R$ contains a red $K_{1,n}$, a contradiction. Thus, we can assume that, any vertex $v \in K_{n+m-1}$ must satisfy $d_R(v) = n-1$ and $d_B(v) = m-1$. Choose the point outside of $K_{n+m-1}$. In order to avoid a red $K_{1,n}$, this vertex cannot be adjacent in red to any vertex of $K_{n+m-1}$. Furthermore, if this vertex $v_0$ is adjacent to some vertex in blue, then since $n\leq m-2$, $\Gamma_B(v_0)$ will contain a red complete graph of order at least $n+1$, a contradiction. Therefore, if the vertex outside of $K_{n+m-1}$ is adjacent in any colour to a vertex of $K_{n+m-1}$, we will get a  red $K_{1,n}$ or a blue $K_{1,m}+e$.  Hence, $r_*(K_{1,n},K_{1,m}+e) \leq 1$. Since by definition,  $r_*(K_{1,n},K_{1,m}+e) \geq 1$, we  conclude that  $r_*(K_{1,n},K_{1,m}+e) = 1$.

\vspace{14pt}
\noindent  \textbf{Case 3.} \textit{$n > m-2$ }

\vspace{8pt}
\noindent Consider the regular standard coloring of $K_{2n}=H_R \oplus H_B$ given in Case 4 of Lemma 1. Extend this coloring to a coloring of $K_{2n} \sqcup K_{1,n}$ such that the new vertex (say $x$) of degree $n$ is adjacent in blue to  all vertices of one partite set of $H_B\cong K_{n,n}$ (see Figure 4). Observe that,  $H_R$ has no $K_{1,n}$. Furthermore, $H_B$ has no $K_{1,m}+e$ since it has no blue $C_3$. That is, $K_{2n} \sqcup K_{1,n} \not \rightarrow (K_{1,n}, K_{1,m}+e)$. Hence, $r_*(K_{1,n},K_{1,m}+e) \geq n+1$.

\vspace{8pt}
\noindent Next we show that, $r_*(K_{1,n},K_{1,m}+e) \leq n+1$. Suppose there exists a red/blue coloring of $K_{2n} \sqcup K_{1,n+1}$ such that $H_R$ contains no $K_{1,n}$ and $H_B$ contains no $K_{1,m}+e$. Let us first restrict our attention to a red/blue coloring of $K_{2n}$. In order to avoid a red $K_{1,n}$, any vertex $v \in K_{2n}$ must satisfy $d_R(v) \leq n-1$ and hence $d_B(v) \geq n$. Next, suppose that there exists a vertex $v_0 \in K_{2n}$ such that $d_R(v_0) \leq n-2$. Then,  $d_B(v_0) \geq n+1 \geq m$. In order to avoid a blue $K_{1,m}+e$, all vertices of $\Gamma_B(v_0)$ must be adjacent to each other in red. Thus, the vertices of $\Gamma_B(v_0)$  will contain a red complete graph of order at least $n+1$. Hence, $H_R$ contains a red $K_{1,n}$, a contradiction. Therefore, we can assume that, any vertex $v \in K_{2n}$ must satisfy $d_R(v) = n-1$ and $d_B(v) = n$. 

\vspace{8pt}
\noindent Let the vertex outside of $K_{2n}$ in $K_{2n} \sqcup K_{1,n+1}$ be denoted by $x$. In order to avoid a red $K_{1,n}$, $x$ cannot be adjacent in red to any vertex of $K_{2n}$. If the vertex $x$ is adjacent to $n+1$ vertices of $K_{2n}$ in blue, then since $n+1\geq m$, $\Gamma_B(x)$ will contain a red complete graph of order at least $n+1$, a contradiction. Hence,  $x$ cannot be adjacent to $n+1$ vertices of $K_{2n}$ in any color.  Therefore, $r_*(K_{1,n},K_{1,m}+e) \leq n+1$. Since by definition,  $r_*(K_{1,n},K_{1,m}+e) \geq n+1$, we can conclude that  $r_*(K_{1,n},K_{1,m}+e) = n+1$.

\begin{center}

\begin{tikzpicture}[line cap=round,line join=round,>=triangle 45,x=1.0cm,y=1.0cm]
\clip(-5.594545454545452,0.51272727272727) rectangle (-0.7763636363636305,7.458181818181823);
\draw [dash pattern=on 5pt off 5pt] (-5.19454545454546,5.003636363636375)-- (-5.230909090909094,1.5672727272727371);
\draw [dash pattern=on 5pt off 5pt] (-5.19454545454546,5.003636363636375)-- (-4.394545454545457,1.5309090909091005);
\draw [dash pattern=on 5pt off 5pt] (-5.19454545454546,5.003636363636375)-- (-3.5763636363636375,1.5127272727272825);
\draw [dash pattern=on 5pt off 5pt] (-4.430909090909093,4.985454545454559)-- (-5.230909090909094,1.5672727272727371);
\draw [dash pattern=on 5pt off 5pt] (-4.430909090909093,4.985454545454559)-- (-4.394545454545457,1.5309090909091005);
\draw [dash pattern=on 5pt off 5pt] (-4.430909090909093,4.985454545454559)-- (-3.5763636363636375,1.5127272727272825);
\draw [dash pattern=on 5pt off 5pt] (-3.667272727272729,4.985454545454559)-- (-5.230909090909094,1.5672727272727371);
\draw [dash pattern=on 5pt off 5pt] (-3.667272727272729,4.985454545454559)-- (-4.394545454545457,1.5309090909091005);
\draw [dash pattern=on 5pt off 5pt] (-3.667272727272729,4.985454545454559)-- (-3.5763636363636375,1.5127272727272825);
\draw [dash pattern=on 5pt off 5pt] (-5.19454545454546,5.003636363636375)-- (-1.2672727272727256,1.476363636363646);
\draw [dash pattern=on 5pt off 5pt] (-4.430909090909093,4.985454545454559)-- (-1.2672727272727256,1.476363636363646);
\draw [dash pattern=on 5pt off 5pt] (-3.667272727272729,4.985454545454559)-- (-1.2672727272727256,1.476363636363646);
\draw [dash pattern=on 5pt off 5pt] (-1.2672727272727256,1.476363636363646)-- (-1.2854545454545439,5.003636363636377);
\draw [dash pattern=on 5pt off 5pt] (-1.2672727272727256,1.476363636363646)-- (-2.9218181818181823,5.003636363636377);
\draw [dash pattern=on 5pt off 5pt] (-1.2672727272727256,1.476363636363646)-- (-2.085454545454545,5.003636363636377);
\draw [dash pattern=on 5pt off 5pt] (-2.085454545454545,5.003636363636377)-- (-2.758181818181822,1.4945454545454546);
\draw [dash pattern=on 5pt off 5pt] (-2.9218181818181823,5.003636363636377)-- (-2.758181818181822,1.4945454545454546);
\draw [dash pattern=on 5pt off 5pt] (-2.9218181818181823,5.003636363636377)-- (-1.94,1.4945454545454622);
\draw [dash pattern=on 5pt off 5pt] (-2.085454545454545,5.003636363636377)-- (-1.94,1.4945454545454622);
\draw [dash pattern=on 5pt off 5pt] (-2.9218181818181823,5.003636363636377)-- (-5.230909090909094,1.5672727272727371);
\draw [dash pattern=on 5pt off 5pt] (-2.085454545454545,5.003636363636377)-- (-5.230909090909094,1.5672727272727371);
\draw [dash pattern=on 5pt off 5pt] (-1.94,1.4945454545454622)-- (-1.2854545454545439,5.003636363636377);
\draw [dash pattern=on 5pt off 5pt] (-1.94,1.4945454545454622)-- (-3.667272727272729,4.985454545454559);
\draw [dash pattern=on 5pt off 5pt] (-1.94,1.4945454545454622)-- (-4.430909090909093,4.985454545454559);
\draw [dash pattern=on 5pt off 5pt] (-1.94,1.4945454545454622)-- (-5.19454545454546,5.003636363636375);
\draw [dash pattern=on 5pt off 5pt] (-2.758181818181822,1.4945454545454546)-- (-1.2854545454545439,5.003636363636377);
\draw [dash pattern=on 5pt off 5pt] (-2.758181818181822,1.4945454545454546)-- (-3.667272727272729,4.985454545454559);
\draw [dash pattern=on 5pt off 5pt] (-2.758181818181822,1.4945454545454546)-- (-4.430909090909093,4.985454545454559);
\draw [dash pattern=on 5pt off 5pt] (-2.758181818181822,1.4945454545454546)-- (-5.19454545454546,5.003636363636375);
\draw [dash pattern=on 5pt off 5pt] (-3.5763636363636375,1.5127272727272825)-- (-2.9218181818181823,5.003636363636377);
\draw [dash pattern=on 5pt off 5pt] (-3.5763636363636375,1.5127272727272825)-- (-2.085454545454545,5.003636363636377);
\draw [dash pattern=on 5pt off 5pt] (-3.5763636363636375,1.5127272727272825)-- (-1.2854545454545439,5.003636363636377);
\draw [dash pattern=on 5pt off 5pt] (-4.394545454545457,1.5309090909091005)-- (-2.9218181818181823,5.003636363636377);
\draw [dash pattern=on 5pt off 5pt] (-4.394545454545457,1.5309090909091005)-- (-2.085454545454545,5.003636363636377);
\draw [dash pattern=on 5pt off 5pt] (-4.394545454545457,1.5309090909091005)-- (-1.2854545454545439,5.003636363636377);
\draw [dash pattern=on 5pt off 5pt] (-5.230909090909094,1.5672727272727371)-- (-1.2854545454545439,5.003636363636377);
\draw (-5.19454545454546,5.003636363636375)-- (-4.430909090909093,4.985454545454559);
\draw (-4.430909090909093,4.985454545454559)-- (-1.2854545454545439,5.003636363636377);
\draw [shift={(-4.440683131897008,4.173526011560707)}] plot[domain=0.809689847466171:2.3080944070070037,variable=\t]({1.0*1.1213347395546562*cos(\t r)+-0.0*1.1213347395546562*sin(\t r)},{0.0*1.1213347395546562*cos(\t r)+1.0*1.1213347395546562*sin(\t r)});
\draw [shift={(-2.8664153322405754,4.129042995839124)}] plot[domain=0.8418946175809634:2.3226855294388913,variable=\t]({1.0*1.1725243321536643*cos(\t r)+-0.0*1.1725243321536643*sin(\t r)},{0.0*1.1725243321536643*cos(\t r)+1.0*1.1725243321536643*sin(\t r)});
\draw [shift={(-4.0581818181818194,3.6712252964427012)}] plot[domain=0.8646426914377155:2.27694996215208,variable=\t]({1.0*1.7511829047902965*cos(\t r)+-0.0*1.7511829047902965*sin(\t r)},{0.0*1.7511829047902965*cos(\t r)+1.0*1.7511829047902965*sin(\t r)});
\draw [shift={(-3.6400000000000015,3.289862258953178)}] plot[domain=0.834078486342882:2.3075141672469126,variable=\t]({1.0*2.3137919638832356*cos(\t r)+-0.0*2.3137919638832356*sin(\t r)},{0.0*2.3137919638832356*cos(\t r)+1.0*2.3137919638832356*sin(\t r)});
\draw [shift={(-3.2399999999999998,2.685454545454556)}] plot[domain=0.8702999568471705:2.2712926967426244,variable=\t]({1.0*3.0321963782105343*cos(\t r)+-0.0*3.0321963782105343*sin(\t r)},{0.0*3.0321963782105343*cos(\t r)+1.0*3.0321963782105343*sin(\t r)});
\draw [shift={(-2.103636363636363,3.9133333333333478)}] plot[domain=0.9270284564927336:2.2145641970970598,variable=\t]({1.0*1.363151563653608*cos(\t r)+-0.0*1.363151563653608*sin(\t r)},{0.0*1.363151563653608*cos(\t r)+1.0*1.363151563653608*sin(\t r)});
\draw [shift={(-2.4665229615745083,3.705417057169648)}] plot[domain=0.8326148768520374:2.3242446557727416,variable=\t]({1.0*1.7550771980853368*cos(\t r)+-0.0*1.7550771980853368*sin(\t r)},{0.0*1.7550771980853368*cos(\t r)+1.0*1.7550771980853368*sin(\t r)});
\draw [shift={(-2.843124090423984,2.389558552440098)}] plot[domain=1.0334137613897836:2.1197394570873334,variable=\t]({1.0*3.04298163226695*cos(\t r)+-0.0*3.04298163226695*sin(\t r)},{0.0*3.04298163226695*cos(\t r)+1.0*3.04298163226695*sin(\t r)});
\draw [shift={(-3.655254237288137,3.242465331278902)}] plot[domain=1.1761915888027186:1.989496284495447,variable=\t]({1.0*1.907787161262603*cos(\t r)+-0.0*1.907787161262603*sin(\t r)},{0.0*1.907787161262603*cos(\t r)+1.0*1.907787161262603*sin(\t r)});
\draw [shift={(-3.248947745168219,3.803350035790994)}] plot[domain=0.8009622871697555:2.356133931844451,variable=\t]({1.0*1.6716470005932818*cos(\t r)+-0.0*1.6716470005932818*sin(\t r)},{0.0*1.6716470005932818*cos(\t r)+1.0*1.6716470005932818*sin(\t r)});
\draw (-5.230909090909094,1.5672727272727371)-- (-1.2672727272727256,1.476363636363646);
\draw [shift={(-3.201505308525176,3.596550366484154)}] plot[domain=3.9269597371464364:5.4519547050299915,variable=\t]({1.0*2.8699211572720276*cos(\t r)+-0.0*2.8699211572720276*sin(\t r)},{0.0*2.8699211572720276*cos(\t r)+1.0*2.8699211572720276*sin(\t r)});
\draw [shift={(-3.9415989684074746,3.33108961960029)}] plot[domain=4.081173965348956:5.284797419012405,variable=\t]({1.0*2.1848044813312626*cos(\t r)+-0.0*2.1848044813312626*sin(\t r)},{0.0*2.1848044813312626*cos(\t r)+1.0*2.1848044813312626*sin(\t r)});
\draw [shift={(-2.0015006615006676,2.4060366660367167)}] plot[domain=4.019527966854847:5.380860959018042,variable=\t]({1.0*1.18464450427212*cos(\t r)+-0.0*1.18464450427212*sin(\t r)},{0.0*1.18464450427212*cos(\t r)+1.0*1.18464450427212*sin(\t r)});
\draw [shift={(-2.783243054097558,4.2364891408309555)}] plot[domain=4.175249545635541:5.214648230982434,variable=\t]({1.0*3.149040938480208*cos(\t r)+-0.0*3.149040938480208*sin(\t r)},{0.0*3.149040938480208*cos(\t r)+1.0*3.149040938480208*sin(\t r)});
\draw [shift={(-3.5325783893525733,3.4830633882247453)}] plot[domain=4.296587112646589:5.083753717469341,variable=\t]({1.0*2.133985384674687*cos(\t r)+-0.0*2.133985384674687*sin(\t r)},{0.0*2.133985384674687*cos(\t r)+1.0*2.133985384674687*sin(\t r)});
\draw [shift={(-3.1333591231463576,3.8018955512571866)}] plot[domain=4.205449584213515:5.18970091433206,variable=\t]({1.0*2.5976855979106186*cos(\t r)+-0.0*2.5976855979106186*sin(\t r)},{0.0*2.5976855979106186*cos(\t r)+1.0*2.5976855979106186*sin(\t r)});
\draw [shift={(-4.357758799584965,2.931619442892501)}] plot[domain=4.143092224018955:5.215775541521828,variable=\t]({1.0*1.6198251114431605*cos(\t r)+-0.0*1.6198251114431605*sin(\t r)},{0.0*1.6198251114431605*cos(\t r)+1.0*1.6198251114431605*sin(\t r)});
\draw [shift={(-3.529885549455663,4.045406159858526)}] plot[domain=4.110839983430501:5.269746275571187,variable=\t]({1.0*3.0057655261643745*cos(\t r)+-0.0*3.0057655261643745*sin(\t r)},{0.0*3.0057655261643745*cos(\t r)+1.0*3.0057655261643745*sin(\t r)});
\draw [shift={(-2.740095693779903,3.131387559808631)}] plot[domain=4.235516555731295:5.167040097242866,variable=\t]({1.0*1.8219234881760282*cos(\t r)+-0.0*1.8219234881760282*sin(\t r)},{0.0*1.8219234881760282*cos(\t r)+1.0*1.8219234881760282*sin(\t r)});
\draw [shift={(-2.382426302453494,3.9959297942030925)}] plot[domain=4.2642145772385955:5.129069923831046,variable=\t]({1.0*2.755318696624454*cos(\t r)+-0.0*2.755318696624454*sin(\t r)},{0.0*2.755318696624454*cos(\t r)+1.0*2.755318696624454*sin(\t r)});
\draw [dash pattern=on 5pt off 5pt] (-3.3763636363636333,7.021818181818189)-- (-5.19454545454546,5.003636363636375);
\draw [dash pattern=on 5pt off 5pt] (-3.3763636363636333,7.021818181818189)-- (-4.430909090909093,4.985454545454559);
\draw [dash pattern=on 5pt off 5pt] (-3.3763636363636333,7.021818181818189)-- (-3.667272727272729,4.985454545454559);
\draw [dash pattern=on 5pt off 5pt] (-3.3763636363636333,7.021818181818189)-- (-2.9218181818181823,5.003636363636377);
\draw [dash pattern=on 5pt off 5pt] (-3.3763636363636333,7.021818181818189)-- (-2.085454545454545,5.003636363636377);
\draw [dash pattern=on 5pt off 5pt] (-3.3763636363636333,7.021818181818189)-- (-1.2854545454545439,5.003636363636377);
\draw (-3.1036363636363595,7.367272727272732) node[anchor=north west] {$x$};
\begin{scriptsize}
\draw [fill=black] (-5.19454545454546,5.003636363636375) circle (1.5pt);
\draw [fill=black] (-4.430909090909093,4.985454545454559) circle (1.5pt);
\draw [fill=black] (-3.667272727272729,4.985454545454559) circle (1.5pt);
\draw [fill=black] (-2.085454545454545,5.003636363636377) circle (1.5pt);
\draw [fill=black] (-1.2854545454545439,5.003636363636377) circle (1.5pt);
\draw [fill=black] (-5.230909090909094,1.5672727272727371) circle (1.5pt);
\draw [fill=black] (-4.394545454545457,1.5309090909091005) circle (1.5pt);
\draw [fill=black] (-3.5763636363636375,1.5127272727272825) circle (1.5pt);
\draw [fill=black] (-1.94,1.4945454545454622) circle (1.5pt);
\draw [fill=black] (-1.2672727272727256,1.476363636363646) circle (1.5pt);
\draw [fill=black] (-2.9218181818181823,5.003636363636377) circle (1.5pt);
\draw [fill=black] (-2.758181818181822,1.4945454545454546) circle (1.5pt);
\draw [fill=black] (-3.3763636363636333,7.021818181818189) circle (1.5pt);
\end{scriptsize}
\end{tikzpicture}
\end{center}
\begin{center}
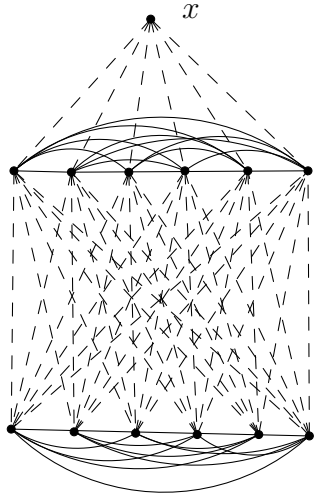
Figure 5. The blue graph of $K_{2n} \sqcup K_{1,n}$ considered in proving $r_*(K_{1,n},K_{1,m}+e) \geq n+1$ when $n=6$ and $m\leq 7$ 
\end{center}

\section{Results and discussion}

In this paper, we proved that the Ramsey number $r(K_{1,n},K_{1,m}+e)$ is $2n+1$ for $n > m-2$. When $n \leq m-2$, $r(K_{1,n},K_{1,m}+e)$ is $n+m+1$ or $n+m$ depending on whether $n$ and $m$ are both even or at least one of them is odd, respectively. 
Furthermore, we showed that the Star critical Ramsey number $r_*(K_{1,n},K_{1,m}+e)$ is $n+1$ for $n > m-2$. When $n \leq m-2$, $r_*(K_{1,n},K_{1,m}+e)$ is $n+m-2$ or 1 depending on whether $n$ and $m$ are both even or at least one of them is odd, respectively. This result is consistent with the known result that, Star critical Ramsey number $r_*(G,H)$ for any two simple graphs $G$ and $H$, satisfies $1 \leq r_*(G,H) \leq r(G,H)-1$. These findings are in agreement with the known result that, Star-critical Ramsey number $r_*(G,H)$ for any two simple graphs $G$ and $H$, satisfies $1 \leq r_*(G,H) \leq r(G,H)-1$.

\vspace{8pt}
\noindent It is interesting to note that when $G=K_{1,n}$ and $H=K_{1,m}+e$, the Star-critical Ramsey number $r_*(G,H)$ achieves the upper bound when \textit{$n$ and $m$ are even and $n \leq m-2$} and  the lower bound when \textit{$n > m-2$}.

\end{document}